\documentclass{amsbook}

\usepackage{amsmath,amstext,amsopn,amsthm,mathrsfs,amssymb}

\usepackage{mathptmx}       
\usepackage{helvet}         
\usepackage{courier}        
\usepackage{type1cm}        
\newtheorem{theo}{Theorem}[chapter]
\newtheorem{lem}{Lemma}[chapter]
\newtheorem{cor}{Corollary}[chapter]
\newtheorem{defi}{Definition}[chapter]

\numberwithin{equation}{chapter}

\allowdisplaybreaks
\begin{document}
\pagenumbering{roman}
\author{MOHAMMAD VALI \hspace{.1cm}SIADAT}
\title{NORM INEQUALITIES FOR INTEGRAL OPERATORS ON CONES\\
THESIS \\Submitted in partial fulfillment of the requirements for the degree of
Doctor of Philosophy in Mathematics in the Graduate College of the
University of Illinois at Chicago,\\
Chicago, Illinois\\
1990}
\maketitle 

\newpage       
\begin{center}
DEDICATIONS
\begin{verse}
I dedicate this thesis to my daughters, Marjan and Banafsheh Siadat who so
 wonderfully supported me during the entire period of its preparation. I am
proud and fortunate to have children like them.
\end{verse}
\end{center}
\bigskip
\begin{center}
\vfill
ACKNOWLEDGEMENTS
\begin{verse}
I would like to thank my thesis advisor, Professor Yoram Sagher, for
introducing me to analysis, for his encouragement and guidance and for all
the help given to me so I could develop mathematically. Without his generous
and unwavering support, this work would not have been possible.\\

I would also like to thank the members of my committee, Professors Calixto
Calder\'{o}n, Nasrollah Etemadi, Jeff Lewis, Charles Lin, and Louis Pennisi
from the University of Illinois at Chicago, and Professor Robert Fefferman
from the University of Chicago for their support and
valuable advice.\\

I was priviledged to have worked with a unique group of talented classmates
who inspired and assisted me with my study and research. Special thanks go
to Wei Cao, Elizabeth Kochneff, Ruby Tan and Kecheng Zhou whose friendship I
will never forget.\\

Finally, I thank Ms. Shirley Roper whose efforts to put this dissertation in
its present form, are greatly appreciated.\\

\hspace{9cm}MVS
\end{verse}
\end{center}
\vfill
\newpage

\tableofcontents

\pagenumbering{arabic}
\mainmatter

\chapter*{SUMMARY}

In this dissertation we explore the $[L^{\mathrm{p}},\ L^{q}]$-boundedness
of certain integral operators on weighted spaces on cones in ${\mathbb R}^{n}.$ These
integral operators are of the type $\displaystyle \int_{V}k(x,\ y)f(y)dy$
defined on a homogeneous cone $V$. The results of this dissertation are then
applied to an important class of operators such as Riemann-Liouville's
fractional integral operators, Weyl's fractional integral operators and
Laplace's operators. As special cases of the above, we obtain an ${\mathbb R}^{n}$
-generalization of the celebrated Hardy's inequality on domains of
positivity. We also prove dual results.

\newpage

\chapter{INTRODUCTION}

Recall Hardy's inequality on the half line: If $1\leq p<\infty$ and $%
\gamma<p-1$, then 
\begin{equation*}
\int_{0}^{\infty}\left(\int_{0}^{x}f(y)dy\right)^{p}x^{\gamma-p}dx\leq
c\int_{0}^{\infty}f^{p}(x)x^{\gamma}dx.
\end{equation*}
This inequality was generalized to self-dual cones in ${\mathbb R}^{n}$ in \cite{Ostrogorski}. We
will further generalize this result to more general cones, and consider
results for related operators.

The present work consists of two chapters. In the first chapter, we review
the geometry and characteristics of convex cones in ${\mathbb R}^{n}.$ In the second
chapter, we prove norm inequalities of some very general operators and their
application to special cases. Here we find conditions for weighted $[L^{p},\
L^{q}]$-boundedness $(1\leq p\leq q\leq \infty )$ of a general $K$-operator
with a special weight $\triangle _{V}^{\delta }(x)$ on a general homogeneous
cone. As a result we obtain a set of sufficient conditions to ensure that
integral operators such as Riemann-Liouville's, Weyl's, Laplace's, Hardy's
and their duals are bounded on such weighted spaces. 

\chapter{ GEOMETRY AND CHARACTERISTICS OF CONVEX CONES}

In this chapter we review some theorems on the geometry and characteristics
of convex cones in ${\mathbb R}^{n}.$ There is no agreement in the literature on the
definition of cones. In some references, see
e.g., \cite{Berman}, \cite{Rockafellar}, cones are closed sets. In others, see e.g., \cite{Vinberg}, cones are
open sets. Because of this lack of uniformity, we present a detailed
exposition of the theory of open cones. Throughout, $c$ stands for a generic
constant which typically depends on the cone under discussion but which may
be different in each occurrence.

\section{Cones}

\begin{defi} A non-empty open subset $V\subset {\mathbb R}^{n}$ is called a cone if
the following conditions hold:
\begin{enumerate}
\item[(i)] if $x\in V$ and $\lambda>0$, then $\lambda x\in V$.
\item[(ii)] if $x, y\in V$, then $x+y\in V$.
\item[(iii)] if $x\neq 0$ and $x\in\overline{V}$, then $-x\not\in\overline{V}$
where $\overline{V}$ denotes the closure of $V.$ Clearly cones are convex
sets.

In the following, we will show that condition (iii) is equivalent to:
\item[(iii)'] $V$ does not contain a line.
\end{enumerate}
\end{defi}
First we introduce some lemmas.

\begin{lem} (see \cite{Bourbaki}). Let $X$ be a topological vector space, and $A\subset X$
be a convex set. Let $x_{0}\in A^{0}$, the interior of $A$, and $x\in%
\overline{A}$. Then every point $y$ in the open segment connecting $x$ and $x_{0}$ lies in $A^{0}$
\end{lem}

\begin{proof}
Let $y$ be a linear combination: $y=(1-\theta)x_{0}+\theta x \; \;
(0<\theta<1).$ Let $\displaystyle \lambda=-\frac{1-\theta}{\theta}.$ Define
a homothety centered at $y$ with ratio $\lambda<0$ as follows: For any $u$
around $y, f(u)$ is defined as $f(u)-y=\lambda(u-y)$. Obviously, $f$ is a
linear transform over $X$ and thus, is one to one and bicontinuous.
Moreover, $f(x_{0})=x.$

Now let $V\subset A$ be a neighborhood of $x_{0}$. Then $f(V)$ is a
neighborhood of $x$. Since $x\in\overline{A}$, there exists $f(z)\in
f(V)\cap A$, where $z\in V.$

Define a homothety centered at $f(z)$ with ratio $\displaystyle \frac{\lambda%
}{\lambda-1}=1-\theta$ as: 
\begin{equation*}
g(v)-f(z)=\frac{\lambda}{\lambda-1}(v-f(z))\ .
\end{equation*}
From $f(z)-y=\lambda(z-y)=\lambda(z-f(z))+\lambda(f(z)-y)$ , we obtain
 $$%
y-f(z)=\displaystyle \frac{\lambda}{\lambda-1}(z-f(z)).$$ This shows that $%
g(z)=y.$

Since $A$ is convex, we have if $v\in A$, then $g(v)\in A$. It follows that $%
y\in A$ and $g(V)\subset A$. Hence $y\in A^{0}$. \end{proof}

\begin{theo} (see \cite{Bourbaki}). Let $X$ be a topological vector space. If $A\subset
X $ is open and convex, then $(\overline{A})^{0}=A.$
\end{theo}

\begin{proof}
 Obviously $A=A^{0}\subset(\overline{A})^{0}.$ So, we only need to show 
$A\supset(\overline{A})^{0}.$ Let $x\in(\overline{A})^{0}.$ We may assume $x=0$%
. Let $V$ be a symmetric neighborhood of $x$ contained in $\overline{A}$.
Since $V\subset\overline{A}$, there exists $y\in V\cap A$. Since $V$ is
symmetric, $-y\in V\subset\overline{A}$. Applying Lemma 2.1, we have that $0$
is an interior point of $A$. That is, $0\in A^{0}=A. $ \end{proof}
The next
lemma characterizes $\overline{V}.$

\begin{lem} Suppose $V$ is an open, convex cone in ${\mathbb R}^{n}.$ Then 
\begin{equation*}
x\in\overline{V}\Leftrightarrow x+V\subset V.
\end{equation*}
\end{lem}
\begin{proof}  Let $x\in\overline{V}$. We have $x+\overline{V}\subset\overline{V}$,
implying $x+V\subset\overline{V}$. Now, since $V$ is open and convex, then
by Theorem 2.1, $(\overline{V})^{0}=V$. Hence, 
\begin{equation*}
x+V=(x+V)^{0}\subset(\overline{V})^{0}=V.
\end{equation*}
Conversely, suppose $x+V\subset V$. If $x\not\in\overline{V}$, then there
exists $\epsilon>0$ such that $B(x,\ \epsilon)\cap V\neq\emptyset.$ Now, let $%
h\in V$ be such that $|h|<\epsilon$; then $x+h\in V$ and $x+h\in B(x,\
\epsilon), \mathrm{a}$ contradiction. 
\end{proof}
We now prove the equivalence of conditions (iii) and (iii)'.

\begin{theo} Let $V$ be a nonempty open subset of ${\mathbb R}^{n}$ which satisfies
(i) if $x\in V$ and $\lambda>0$, then $\lambda x\in V$ and (ii) if $x, y\in
V $, then $x+y\in V$. Then $V$ does not contain a line if and only if $x\neq
0$ and $x\in\overline{V}$ implies $-x\not\in\overline{V}.$
\end{theo}

\begin{proof} Suppose $V$ does not contain a line. If $x\neq 0, \; x \in \bar{V}$ and $-x \in \bar{V},$ then $\lambda x\in\overline{V}$ for all
real $\lambda$. Now if we let $h\in V$, we have by Lemma 2.2 that $\lambda
x+h\in V$. So, $V$ contains a line which is a contradiction.

Conversely, suppose that if $x\in \overline{V}$ and $x\neq 0$, then $%
-x\not\in \overline{V}$. Let $h+\lambda x$ be any line with $h+\lambda x\in V
$ for all $\lambda $. Taking $\lambda =0$, we have $h\in V$. Now, $\forall
\rho >0,\rho (h+\lambda x)\in V$. So, $\rho h+x\in V$ and $\rho h-x\in V$.
Letting $\rho \rightarrow 0^{+}$, we see that $x\in \overline{V}$ and $-x\in 
\overline{V}$. Therefore, $V$ cannot contain a line. 
\end{proof}

\section{ Dual of a Cone}

In this section, we introduce the dual $V^{*}$ of a cone $V$ and discuss the
relationships between $V, V^{*}$ and $(V^{*})^{*}$.

\begin{defi} We define the dual $V^{*}$ of $V$:

$$V^{*}=\{x\in {\mathbb R}^{n} : x\cdot y>0, \;\text{for every}\; y\in\overline{V}, y\neq 0\}.$$

A cone is called self-dual if $V=V^{*}$
\end{defi}

We will show that $V^{*}$ is a non-empty, open convex cone in ${\mathbb R}^{n}$

\begin{lem} Let $(X,T)$ be a topological vector space. Suppose $A$ is a
compact set and $B$ is a closed set with $A\cap B=\emptyset.$
Then there exists an open neighborhood $V$ of $0$ such that $V$ is convex
and 
\begin{equation*}
(A+V)\cap(B+V)=\emptyset.
\end{equation*}
\end{lem}
\begin{proof} Take any $x\in A$. Then $x\in B^{\mathrm{c}}$ where $B^{c}$ is the
complement of $B$. So, $0\in-x+B^{\mathrm{c}}$. There exists a balanced
neighborhood $V_{x}$ of $0$ such that $V_{x}+V_{x}+V_{x}\subset-x+B^{\mathrm{%
c}}.$

Therefore, $x+V_{x}+V_{x}\subset(V_{x}+B)^{c}.$ Now
$x+V_{x}$ is open and $\displaystyle \bigcup_{x\in A}(x+V_{x})$ covers $
A. $ So, there is a finite subcover, i.e., there exists $x_{1},\cdots,  x_{m}\in A$ such that 
\begin{equation*}
\bigcup_{j=1}^{m}\{x_{j}+V_{xj}\}\supset A.
\end{equation*}

Let $V=\cap_{j=1}^m V_{x_j}$ be an open neighborhood of $0$. We have 
\begin{eqnarray*}
A+V&\subset&\bigcup_{j=1}^{m}(x_{j}+V_{x_j}+V)\subset\bigcup_{j=1}^{m}(x_{j}+V_{x_j}+V_{x_j})\\
&\subset&\bigcup_{j=1}^{m}(V_{x_j}+B)^{c}
=[\bigcap_{j=1}^{m}(V_{x_j}+B)]^{\mathrm{c}}
\subset(V+B)^{\mathrm{c}}.
\end{eqnarray*}
Hence, 
\begin{equation*}
(A+V)\cap(B+V)=\emptyset.
\end{equation*}
\end{proof}
\begin{lem} (Separation Theorem) (see \cite{Berman} and \cite{Rockafellar}).

Let $(X,\ T)$ be a topological vector space. Let $A, B$ be non-empty convex
sets such that $A\cap B=\emptyset$. Then we have the following:
\begin{enumerate}
\item[(1)] If $A$ is open, then there exists a continuous linear functional $L$
such that 
\begin{equation*}
\mathrm{Re} L(a)<\inf\{\mathrm{Re} L(b):b\in B\},\;\forall a\in A.
\end{equation*}
\item[(2)] If $X$ is locally convex, $A$ is compact and $B$ is closed, then there
exists a continuous linear functional $L$ such that 
\begin{equation*}
\sup\{\mathrm{Re} L(a)\ :\ a\in A\}<\inf\{\mathrm{Re} L(b):b\in B\}.
\end{equation*}
\end{enumerate}
\end{lem}
\begin{proof}

 (1) Let us consider first the case of $(X,\ T)$ being a topological
vector space over R. Choose $a_{0}\in A, b_{0}\in B$. Let $x_{0}=b_{0}-a_{0}$
and $C=x_{0}+A-B$. Clearly, $0\in C$. Also, since $A-b$ is open, $C=x_{0}+%
\displaystyle \bigcup_{b\in B}(A-b), C$ is open.

Now let $x_{1}, x_{2}\in C$ and $\alpha, \beta\geq 0$ with $\alpha+\beta=1$.
We have $x_{1}=x_{0}+a_{1}-b_{1}, x_{2}=x_{0}+a_{2}-b_{2}$ for some $a_{1},
a_{2}\in A$ and $b_{1}, b_{2}\in B$. So, $\alpha x_{1}+\beta x_{2}=(\alpha
x_{0}+\beta x_{0})+ (\alpha a_{1}+\beta a_{2})-(\alpha b_{1}+\beta b_{2})\in
C$. Therefore, $C$ is convex.

Now let $p(x)=\displaystyle \inf\{t>0\ :\ x\in tC\}$ be the Minkowski's
functional. We have that $p(x)<1\; \forall x\in C$. Also, since $A\cap B=\emptyset,
x_{0}\not\in C$. Since $x_{0}\not\in C$, for any $t<1, x_{0}\not\in tC$
from convexity of $C$ and $0\in C$. It follows that $p(x_{0})\geq 1$. Now
consider $\{tx_{0}\}_{t\in R}=M$ and define $\ell(tx_{0})=t$, so that $%
\ell(x_{0})=1$. We have $\ell(tx_{0})\leq p(tx_{0})$.

Since $M\subset X$, by Banach-Hahn theorem there exists a linear functional $%
L$ defined on $X$ such that
$$L|_{M}=\ell\;  \text{and} \; L(x)\leq p(x) \; \forall x\in X.$$
Since $p(x)<1, \; \forall x\in C$, we have $L(x)<1 \; \forall x\in C$. Also, since $
L$ is linear, i.e., $L(-x)=-L(x),$ we have that $L(x)>-1 \; \forall x\in-C$.
So, $|L(x)|<1 \; \forall x\in C\cap(-C)$ . Hence, $L$ is continuous.

Now let $a\in A, \,b\in B$. We have since $l(x_{0})=1,\; L(a)-L(b)+1=L(a-b+x_{0})\leq p(a-b+x_{0})<1$.
Therefore, $L(a)<L(b)$ . That is to say, for all $a\in A,$ 
\begin{equation*}
L(a)\leq\inf\{L(b)\ :\ b\in B\}.
\end{equation*}
Since $A$ is open, if $a\in A$ then if $|\lambda-1|<e, \lambda a\in A$ and
thus $L(a)<\displaystyle \inf\{L(b):b\in B\}.$ If $(X,\ T)$ is a topological
vector space over $\mathrm{C}$ we consider it first as a t.v.s. over $%
\mathrm{R}$, and construct $L_{0}$ as above. We then define $%
L(x)=L_{0}(x)-iL_{0}(ix)$ and get the claim.

(2) Since $X$ is locally convex, by Lemma 2.3 there exists an open
neighborhood $V$ of $0$ such that $V$ is convex and 
\begin{equation*}
(A+V)\cap(B+V)=\emptyset.
\end{equation*}
Now, $A+V, B+V$ are convex, so by part (1), there exists a continuous linear
functional $L$ such that 
\begin{equation*}
\mathrm{Re} L(x)<\inf\{\mathrm{Re} L(y)\ :\ y\in B+V\}\leq\inf\{\mathrm{Re}
L(y):y\in B\},
\;\forall x\in A+V.
\end{equation*}
Also, since $A$ is compact, there exists $a_{0}\in A$ so that 
\begin{equation*}
\mathrm{Re} L(a_{0})=\sup_{a\in A}\mathrm{Re} L(a)\ .
\end{equation*}
Hence, $\displaystyle \sup_{a\in A}\mathrm{Re} L(a)=\mathrm{Re}
L(a_{0})<\inf\{\mathrm{Re} L(y)\ :\ y\in B\}.$
\end{proof}

\begin{cor} (see \cite{Rockafellar}). Let $V$ be a cone and let $a\in {\mathbb R}^{n}, \;
a\not\in\overline{V}$. Then, there exists a hyperplane separating $a$ and $%
\overline{V}.$
\end{cor}

\begin{proof} By part (2) of Lemma 2.4, there exists a fixed point $x_{0}\in {\mathbb R}^{n}$
such that 
\begin{equation*}
\sup\{(x\cdot x_{0})\ :\ x\in\overline{V}\}<(x_{0}\ .\ a)\ .
\end{equation*}
Now let $\displaystyle \alpha=\sup\{(x\cdot x_{0})\ :\ x\in\overline{V}\}$.
Since $\overline{V}$ is closed, $0\in\overline{V}$, so that $\alpha\geq 0$.
Now, for any $\beta>0,$ 
\begin{equation*}
\beta(x\cdot x_{0})=(\beta x\cdot x_{0})\leq\alpha.
\end{equation*}
This in turn implies that $(x\cdot x_{0})\leq 0,  \;\forall x\in\overline{V}$.
So, $\alpha\leq 0$. Therefore, $\alpha=0$. It follows that $(x\cdot
x_{0})\leq 0,\; \forall x\in\overline{V}$ and $(x_{0}\cdot a)>0.$
\end{proof}

\begin{defi} Let $V$ be a cone. Let $\Sigma_{\overline{V}}$ denote the
intersection of $\overline{V}$ and the unit sphere, i.e., $\Sigma_{\overline{%
V}}=\{x\in\overline{V}\ :\ |x|=1\}$. The convex hull of $\Sigma_{\overline{V}%
}$, denoted $CVX(\Sigma_{\overline{V}})$), is defined a\textit{s} the
intersection of all convex sets containing $\Sigma_{\overline{V}}.$
\end{defi}

\begin{theo} (see \cite{Taylor}). The convex hull of a set $S$ consists of all points
which are expressible in the form $\displaystyle \sum_{i=1}^N\alpha_{i}x_{i}$
where $x_i\in S, \alpha_{i}>0$ and $\displaystyle \sum_{i=1}^{N}\alpha_{i}=1.$
\end{theo}

\begin{proof} Trivial.
\end{proof}

\begin{theo} (Carath\'{e}odory's Theorem) (see \cite{Rockafellar}). Let $S$ be any set of
points and directions in ${\mathbb R}^{n}$, and let $C=CVX(S)$. Then $x\in C$ if and
only if $x$ can be expressed as a convex combination of $n+1$ of the points
and directions in $S$ (not necessarily distinct). In fact $C$ is the union
of all the generalized $d$-dimensional simplices whose vertices belong to $S$, where $d=\dim C.$
\end{theo}

\begin{proof} Let $S_{0}$ be the set of points in $S$ and $S_{1}$ the set of
directions in $S$. Let $S'_{1}$ be a set of vectors in ${\mathbb R}^{n}$ such
that the set of directions of the vectors in \textit{S}\'{\i} is $S_{1}$.
Let $S^{\prime }$ be the subset of ${\mathbb R}^{n+1}$ consisting of all the vectors
of the form (1, \textit{x}) with $x\in S_{0}$ or of the form $(0,\ x)$ with $%
x\in S\text{\'{\i}}.$ Let $K$ be the convex cone generated by $S^{\prime}.$
$CVX(S)$ can be identified with the intersection of $K$ and the hyperplane $%
\{(1,\ x)|x\in {\mathbb R}^{n}\}$. Translating the statement of the theorem into this
context in ${\mathbb R}^{n+1}$, we see that it is only necessary to show that any
non-zero vector $y\in K$, which is in any case a non-negative linear
combination of elements of $S^{\prime }$, can actually be expressed as a
non-negative linear combination of $d+1$ linearly independent vectors of $%
S^{\prime }$, where $d+1$ is the dimension of $K(=$ the dimension of
the subspace of ${\mathbb R}^{n+1}$ generated by $S^{\prime }$). The argument is
algebraic, and it does not depend on the relationshhip between $S^{\prime }$
and $S$. Given $y\in K$, let $y_{1}$, . . . , $y_{m}$ be vectors in $%
S^{\prime }$ such that $y=\lambda_{1}y_{1}+\ldots+\lambda_{m}y_{m}$, where
the coefficients $\lambda_{i}$ are all non-negative. Assuming the vectors $%
y_{i}$ are not themselves linearly independent, we can find scalars $\mu_{1}$%
, . .. , $\mu_{m}$, at least one of which is positive, such that $%
\mu_{1}y_{1}+\ldots+\mu_{m}y_{m}=0$. Let $\lambda$ be the greatest scalar
such that $\lambda\mu_{i}\leq\lambda_{i}$ for $i=1$, . . . , $m$, and let $\lambda_i = \lambda_i - \lambda \mu_i.$ Then 
\begin{equation*}
\lambda_{1}^{\prime }y_{1}+.\text{ . . }+\lambda_{m}^{\prime
}y_{m}=\lambda_{1}y_{1}+.\text{ . . }+\lambda_{m}y_{m}-\lambda\
(\mu_{1}y_{1}+.\ .\ .\ +\mu_{m}y_{m})=y.
\end{equation*}
By the choice of $\lambda$, the new coefficients $\lambda'_i$
are non-negative, and at least one of them is 0. We therefore have an
expression of $y$ as a non-negative linear combination of fewer than $m$
elements of $S^{\prime }$. If these remaining elements are not linearly
independent, we can repeat the argument and eliminate another of them. After
a finite number of steps, we get an expression of $y$ as a non-negative
linear combination of linearly independent vectors $z_{1}, \cdots, z_{r}$
of $S^{\prime }$. Then $r\leq d+1$ by definition of $d+1$. Choosing
additional vectors $z_{r+1}, \cdots, z_{d+1}$ from $S^{\prime }$ if
necessary to make a basis for the subspace generated by $S^{\prime }$, we
get the desired expression of $y$ by adding the term $0z_{r+1}+%
\ldots+0z_{d+1}$ to the expression in terms of $z_{1},\cdots, z_{r}.$ 
\end{proof}

\begin{cor} (see \cite{Rockafellar}). If $S$ is a bounded set of points in ${\mathbb R}^{n}$, then 
$\overline{CVX(S)}=CVX(\overline{S})$ . In particular, if $S$ is closed and
bounded, then $CVX(S)$ is closed and bounded.
\end{cor}

\begin{proof} Let $m=(n+1)^{2}$, and let $Q$ be the set of all vectors of the form 
\begin{equation*}
(\lambda_{0},\ .\ .\ .\ \lambda_{n},\ x_{0},\ .\ .\ .\ x_{n})\in {\mathbb R}^{m}
\end{equation*}
such that the components $\lambda_{i}\in R$ and $x_{i}\in {\mathbb R}^{n}$ satisfy 
\begin{equation*}
\lambda_{i}\geq 0, \;  \lambda_{0}+\ldots+\lambda_{n}=1,\ x_i \in
\overline{S}.
\end{equation*}
The image of $Q$ under the continuous mapping 
\begin{equation*}
\theta:\ (\lambda_{0}, \ldots,\ \lambda_{n}, x_{0}, \cdots, x_{n})\rightarrow\lambda_{0}x_{0}+\ldots+\lambda_{n}x_{n}
\end{equation*}
from ${\mathbb R}^{m}$ to ${\mathbb R}^{n}$ is $CVX(\overline{S})$ by Carath\'{e}odory's
Theorem. If $S$ is bounded in ${\mathbb R}^{n}, Q$ is closed and bounded in ${\mathbb R}^{m}$,
and hence the image of $Q$ under $\theta$ is closed and bounded too. Then 
\begin{equation*}
CVX(\overline{S})=\overline{CVX(\overline{S})}\supset\overline{CVS(S)}.
\end{equation*}
Of course, in general 
\begin{equation*}
\overline{CVX(S)}=CVX(\overline{CVX(S)})\supset CVX(\overline{S})
\end{equation*}
so the commutativity of ``convex hull'' and ``closure'' follows.
\end{proof}

\begin{lem} If $V$ is cone, then $0\not\in CVX(\Sigma_{\overline{V}})$ .
\end{lem} 
\begin{proof} Suppose $0\in CVX(\Sigma_{\overline{V}})$ . Then

$0=\displaystyle \sum_{i=1}^{N}\alpha_{i}x_{i}$, where $\alpha i>0, %
\displaystyle \sum_{i=1}^{N}\alpha_{i}=1$ and $x_{i}\in\Sigma_{\overline{V}%
}. $

So, $\displaystyle \alpha_{1}x_{1}=-\sum_{j=2}^{N}\alpha_{j}x_{j}$. Since $%
\alpha_{1}x_{1}\in\overline{V}$ and $\displaystyle \sum_{j=2}^{N}%
\alpha_{j}x_{j}\in\overline{V}$, we get contradiction. 
\end{proof}

\begin{lem} $CVX(\Sigma_{\overline{V}})=\overline{CVX(\Sigma_{V})}.$
\end{lem}

\begin{proof} By Corollary 2.2, it suffices to show that $\overline{\Sigma_{V}}%
=\Sigma_{\overline{V}}$. Clearly $\Sigma_{\overline{V}}$ is closed since it
is the intersection of $\overline{V}$ and the unit ball so that $\overline{%
\Sigma_{V}}\subset\Sigma_{\overline{V}}.$

To show $\Sigma_{\overline{V}}\subset\overline{\Sigma_{V}}$, let $x\in%
\overline{V}, |x|=1$ and show that there exists $x_{n}\in V, |x_{n}|=1$
such that $x_{n}\rightarrow x$. Since $x\in\overline{V}$, there exists $%
y_{n}\in V$ so that $y_{n}\rightarrow x$. Let $x_{n}=\displaystyle \frac{%
y_{n}}{|y_{n}|}$. Then $|x_{n}|=1$. Also $x_{n}\in V$ since $V$ is a cone.
Now $x_{n}=\displaystyle \frac{y_{\mathfrak{n}}}{|y_{n}|}\rightarrow\frac{x}{%
|x|}=x.$
\end{proof}

\begin{lem} If $V$ is a cone, then $-CVX(\Sigma_{\overline{V}})\cap\overline{V%
}=\emptyset.$
\end{lem}

\begin{proof} Clearly $\overline{V}\cap(-\overline{V})=\{0\}$, otherwise $x\in%
\overline{V}, x\in (-\overline{V})$ , i.e., $-x\in\overline{V}$
contradicting condition (iii) of the definition of a cone. Now since 
$$
-CVX(\Sigma_{\overline{V}})\subset-\overline{V},\; -CVX(\Sigma_{\overline{V}})\cap\overline{V}\subset\{0\}.$$
 But, by Lemma 2.5, $0\not\in CVX(\Sigma_{\overline{V}})$.
\end{proof}

\begin{theo} $V^{*}\neq\emptyset.$
\end{theo}

\begin{proof} From Lemma 2.6 we have that $CVX(\Sigma_{\overline{V}})$ is compact.
By Lemma 2.4, there exists $x_{0}\in {\mathbb R}^{n}$ such that 
\begin{equation*}
\sup\{(x_{0}\cdot y):y\in-CVX(\Sigma_{\overline{V}})\}<\inf \{(x_{0}\cdot x):x\in\overline{V}\}=\beta.
\end{equation*}
But $\beta\leq 0$, because $0\in\overline{V}$. Also, $\beta$ cannot be
negative, otherwise it will be $-\infty$. So, $\beta=0$. Hence, $(x_{0}\cdot y)
<0$ for all $y\in-CVX(\Sigma_{\overline{V}}),$ implying that $(x_{0} \cdot y)
>0$ for all $y\in CVX(\Sigma_{\overline{V}})$ . So, $x_{0}\in V^{*}$ and
hence $V^{*}\neq\emptyset. $
\end{proof} 

\begin{theo} $V^{*}$ is open.
\end{theo}

\begin{proof} Let $y\in V^{*}.$ Show that there exists $\delta>0$ such that $%
y+\delta B\subset V^{*}$ where $B$ is an open unit ball, i.e., $%
B=\{y:|y|<1\}.$

Take $y^{\prime }\in B$, and show $y+\delta y^{\prime} \in V^*$, i.e., show that $
(y+\delta y^{\prime })\cdot x>0, \; \forall x\in\overline{V}$. Now consider $%
y\cdot x$ where $x\in\Sigma_{\overline{V}}. y\cdot x$ is continuous in $x$
and thus attains its extremum on the compact set $\Sigma_{\overline{V}}.$
So, $y\cdot x\geq\lambda>0$ for some $\lambda>0$, and all $x\in\Sigma_{%
\overline{V}}$. Now choose any number $\delta, 0<\delta<\lambda$. For $%
x\in\Sigma_{\overline{V}}$ and using Schwarz inequality, $(y+\delta
y^{\prime })\cdot x=y\cdot x+\delta y^{\prime }\cdot x\geq
\lambda-\delta|y^{\prime }\cdot x|\geq\lambda-\delta>0$. For $x\in\overline{V
},\frac{x}{|x|}\in\Sigma_{\overline{V}}, x\neq 0, (y+\displaystyle \delta
y^{\prime })\cdot x=[(y+\delta y^{\prime })\cdot\frac{x}{|x|}]|x|>0.$ So, $%
V^{*}$ is open. 
\end{proof}

\begin{theo} $V^{*}$ is a cone.
\end{theo} 

\begin{proof} Conditions (i) and (ii) of the definition of a cone trivially hold.
To show (iii), suppose $x\in\overline{V}^{*}, \; x\neq 0,$ and $-x\in\overline{V%
}^{*}.$ So, $y\cdot x\geq 0, \; \forall y\in V$ and $y\cdot(-x)\geq 0, \; \forall
y\in V.$ Hence, $y\cdot x=0, \;\forall y\in V$. But, since $V$ is open, there
exists $\eta>0$ such that $y+\eta z\in V$ for all $|z|<1$, implying $(y+\eta
z)\cdot x=0$. Therefore $z\cdot x=0$ which implies that $x=0.$\end{proof}
 We denote by $V^{**}$ the dual $(V^{*})^{*}$ of $V^{*}$.

\begin{theo} $V=V^{**}$ (see \cite{Berman} and \cite{Rockafellar} for the corresponding theorem for
closed cones).
\end{theo}

\begin{proof} Show first that $V\subset V^{**}.$ It suffices to show that $x\cdot
y>0,\; \forall y\in\Sigma_{\overline{V^{*}}}$ and $x\in V.$ Clearly, by
continuity of inner product, for all $x\in V$, if $y\in\Sigma_{\overline{
V^{*}}}, y\neq 0$, then $x\cdot y\geq 0.$

 Let us assume $x \cdot y=0$ for some $
y\in\Sigma_{\overline{V^{*}}}, y\neq 0$. Since $V$ is open, there exists $
\rho>0$ such that $x+\rho z\subset V,  \;\forall |z|=1$. So, $(x+\rho z) y\geq 0,$
 implying $\rho z$ . $y\geq 0$. Also, $(x-\rho z)\cdot y\geq 0$, implying $
\rho z\cdot y\leq 0$. As a result, $z \cdot y=0, \; \forall |z|=1$, and so $y=0$. So, $
x\cdot y>0, \; \forall y\in\Sigma_{\overline{V^{*}}}$ and we have shown that $
V\subset V^{**}$.

To show $V^{**}\subset V$, by Theorem 2.1 it suffices to show that $
V^{**}\subset\overline{V}$. Suppose, $x\in V^{**}$, but $x\not\in\overline{V}
$. Then by Corollary 2.1, there exists $y\in {\mathbb R}^{n}$ such that $y\cdot x<0$,
and $y\cdot u\geq 0, \; \forall u\in\overline{V}.$

We will show that $y\in\overline{V^{*}}$. Take $z\in V^{*}.$ For every $%
\rho>0, (y+\rho z)\cdot u>0\; \forall u\in\overline{V}, u\neq 0$. So, $y+\rho
z\in V^{*}$. By letting $\rho\rightarrow 0^{+}$, we get that $y\in\overline{%
V^{*}}$. But, we had that $y\cdot x<0,$ contradicting the fact that $x\in
V^{**}$. Therefore $x\in\overline{V}$ and so $V=V^{**} .$
\end{proof} 

\begin{cor}  Suppose $V_{1}$ and $V_{2}$ are cones. Then,
\begin{enumerate}
\item[ (i)] $V_{1}\subset
V_{2}\Leftrightarrow V_{2}^{*}\subset V_{1}^{*}$.
\item[(ii)] $V_{1}=V_{2}\Leftrightarrow V_{1}^{*}=V_{2}^{*}$.
\end{enumerate}
\end{cor}

\begin{proof}

(i) Assume $V_{1}\subset V_{2}$ and let $x\in V_{2}^{*}.$ We have
that $x\cdot y>0$ for every $y\in\overline{V}_{2}\supset\overline{V}_{1},$
so that $x\in V_{1}^{*}$. Conversely, if $V_{2}^{*}\subset V_{1}^{*}$, then by the above result $%
(V_{1}^{*})^{*}\subset(V_{2}^{*})^{*}$, i.e., $V_{1}\subset V_{2}.$

(ii) If $V_{1}^{*}=V_{2}^{*}$, then $V_{1}^{**}=V_{2}^{**}$, i.e., $V_{1}=V_{2}.$

If $V_{1}=V_{2}$, then trivially $V_{1}^{*}=V_{2}^{*}.$
\end{proof}

\section{The Characteristic Function of a Cone}

\begin{defi} (see \cite{Koeche}). The characteristic function $\phi_{V}(x)$ of the
cone is defined as 
\begin{equation*}
\phi_{V}(x)=\int_{V^*}\ e^{-x\cdot y}dy,\ x\in V.
\end{equation*}
\end{defi}
\begin{lem}
 If $x\in V$ and $r(x)=$dist \{$x,\ \partial V\}$, where $\partial
V$ denotes the boundary of $V,$ then $x\cdot y\geq r(x), \; \forall y\in\Sigma_{
\overline{V^{*}}}.$
\end{lem}

\begin{proof} For all $r,\;  0<r<r(x), z=x-yr\in V=V^{**}.$

Now, since $y\in\Sigma_{\overline{V^{*}}} \subset\overline{V^{*}}$, it follows
that $0<y z=y x-r$. Hence, $y x>r,$ implying $y\cdot x\geq r(x)$ . 
\end{proof}

\begin{theo} (see \cite{Vinberg}).
$$\displaystyle \phi_{V}(x)=\int_{V^{*}}e^{-x
\cdot y}dy<\infty,$$ 
for every $x\in V.$
\end{theo}

\begin{proof} Using polar coordinates and Lemma 2.8 as well as change of variables, 
$y=\rho t, \;\rho=|y|$, we get, 
\begin{eqnarray*}
\int_{V^{*}}e^{-x\cdot y}dy&=&\int_{\Sigma_{V^{*}}}\int_{0}^{\infty}e^{(-x\cdot t)\rho}\rho^{n-1}d\rho dt\\
&\leq&\left(\int_{\Sigma_{\overline{V^{*}}}} dt\right) \left(\int_{0}^{\infty}e^{-r(x)\rho}\rho^{n-1}dp\right)=\mu(\Sigma_{\overline{V^{*}}})\frac{\Gamma(n)}{r^{n}(x)}<\infty,\\
\end{eqnarray*}
where $\mu(\Sigma_{\overline{V^{*}}})$ is the measure of $\Sigma_{\overline{V^{*}}}$ and $\Gamma(n)$ is the gamma function evaluated at $n.$
\end{proof}

\begin{theo} $\phi_{V}(x)$ is strictly $\log$ convex.
\end{theo}
\begin{proof} For every $0\leq\theta\leq 1$, and $x_{0}, x_{1}\in V$, we have by H\"{o}lder's inequality that 
\begin{eqnarray*}
\phi_{V}((1-\theta)x_{0}+\theta x_{1})&=&\int_{V^{*}}e^{-(\langle
1-\theta)x_{0}+\theta x_{1})\cdot y}dy\\
&\leq&\left(\int_{V^{*}}e^{-x_{0}\cdot y}dy\right)^{1-\theta}\left(\int_{V^{*}}e^{-x_{1}\cdot y}dy\right)^{\theta}.
\end{eqnarray*}
Since $e^{-x_{0}\cdot y}$ and $e^{-x_{1}\cdot y}$ are not proportional, we
get the strict inequality in above and so $\phi_{V}(x)$ is strictly $\log$ convex. 
\end{proof} 

\begin{lem} Suppose $x_{0}\in\partial V^{*}.$ Then there exists $y_{0}\in%
\overline{V}$ such that $x_{0}\cdot y_{0}=0.$
\end{lem} 
\begin{proof} If for every $y\in\Sigma_{\overline{V}}, \; x_{0}\cdot y>0,$
then 
\begin{equation*}
\min\{\ (x_{0}\ .\ y):y\in\Sigma_{\overline{V}}\}=\delta>0
\end{equation*}
and so, $x_{0}\in V^{*}$, a contradiction. 
\end{proof}

\begin{theo} (see \cite{Vinberg}).

$$\lim _{x \rightarrow \partial V \atop x \in V} \phi_{V}(x)=\infty.$$
\end{theo}

\begin{proof} Take $x_{0}\in\partial V$. By Fatou's lemma, 
\begin{equation*}
\liminf _{x \rightarrow x_0 \atop x \in V} \int_{V^{*}} e^{x \cdot y} d y \geq \int_{V^{*}} e^{-x_{0} \cdot y} d y.
\end{equation*}
So, it suffices to show that $\displaystyle \int_{V^*} e^{-x_{0}\cdot
y}dy=\infty$. By Lemma 2.9, since $V^{**}=V$, we can choose $y_{0}\in \overline{V^{*}}, y_{0}\neq 0$, such that $x_{0}\cdot y_{0}=0.$

Now, take a closed ball $B\subseteq V^{*}$ and let $L=B+\{\lambda
y_{0}\}_{\lambda>0}$. By Lemma 2.2, $L\subset V^{*}.$ Now, for every $y\in L,x
y=z+\lambda y_{0}$ where $z\in B$; also since $B$ is closed and bounded $%
e^{-x_{0}\cdot z}$ will have a minimum on $B.$

So, $e^{-x_{0}\cdot y}=e^{-x_{0}\cdot(z+\lambda y_{0})}=e^{-x_{0}\cdot
z}\geq c>0$. Hence
\begin{equation*}
\int_{V^{*}}e^{-x_{0}\cdot y}dy\geq\int_{L}e^{-x_{0}\cdot
y}dy\geq\int_{L}cdy=\infty.
\end{equation*}
\end{proof}

From Theorem 2.11 it is easy to see that if $V_{1}$ and $V_{2}$ are cones
with $V_{1}\cap V_{2}\neq\emptyset$ such that $\phi_{V_{1}}(x)=\phi_{V_{2}}(x)$
on $V_{1}\cap V_{2}$, then $V_{1}=V_{2}$; otherwise there exists $x^{\prime
}\in V_{1}, x^{\prime }\in\partial V_{2}$ and then $\phi_{V_{2}}(x^{\prime
})=\infty$ while $\phi_{V_{1}}(x^{\prime })<\infty.$

\section{Automorphism Group of a Cone and Homogeneous Cones}

Let $V$ be a cone. Since $V$ is an open set, any matrix which maps $V$ onto $%
V$ is regular. The group of all matrices $A$ which map $V$ onto $V$ is
called the automorphism group of $V$ and denoted by $G(V)$ .

\begin{defi} The cone $V$ is said to be homogeneous if $G(V)$ is
transitive, i.e., for every $x, y\in V$, there exists $A\in G(V)$ such that 
$y=Ax.$
\end{defi}

From now on we will consider only homogeneous cones. We will therefore write
``cone'' for ``homogeneous cone''.

\begin{theo}  $G(V)=G(\overline{V})$ . 
\end{theo} 

\begin{proof} 

Let $A\in G(V)$.

If $s\in V$, then $As\in V$ and $A^{-1}s\in V$. By Lemma 2.2, for any $x\in \overline{V}, \;x+As\in V,$ implying 
\begin{equation*}
A^{-1}x+s\in V.
\end{equation*}
Hence, $A^{-1}x\in\overline{V}.$
Also, we have that $x+A^{-1}s\in V$, implying 
\begin{equation*}
Ax+s\in V.
\end{equation*}
Therefore, $Ax\in\overline{V}$ and as a result $A\overline{V}=\overline{V}$
and $A\in G(\overline{V})$.\\

Conversely, let $B\in G(\overline{V}).$ We have that
$$ B:\overline{V}\rightarrow\overline{V} \; \text{and}\; B^{-1}: \overline{V} \rightarrow\overline{V}.$$

Since $B$ and $B^{-1}$ are continuous,

$$B: (\overline{V})^{0}\rightarrow(\overline{V})^{0}\; \text{and}\; B^{-1}: (\overline{V})^{0}\rightarrow(\overline{V})^{0},$$
i.e,  
$$B:V\rightarrow V \; \text{and}\; B^{-1} : V\rightarrow V,$$
and $B\in G(V)$.
\end{proof}

\begin{theo} $A\in G(V)\Leftrightarrow A^{t}\in G(V^{*})$ , where $A^{t}$
is the transpose of $A.$
\end{theo} 

\begin{proof} Let $A\in G(V)$ and $x\in V^{*}$. Since by Theorem 2.12, $G(V)=G(
\overline{V})$), for every $s\in\Sigma_{\overline{V}}, \;A^{\mathrm{t}}x\cdot
s=x\cdot As>0$. Hence, $A^{t}x\in V^{*}.$ Similarly, for every $s\in\Sigma_{%
\overline{V}}, (A^{t})^{-1}x\cdot s= x\cdot A^{-1}s>0$ and thus $%
(A^{t})^{-1}x\in V^{*}.$ Therefore, $A^{t}\in G(V^{*})$ .

Conversely if $A^{t}\in G(V^{*}),$ then $A=(A^{t})^{t}\in G(V^{**})=G(V)$ . 
\end{proof}

\section{Domains of Positiviity}

In this section we shall discuss domains of positivity. Domains of
positivity serve as a generalization of self-dual cones.

\begin{defi} (see \cite{Rothaus}). $D\subset {\mathbb R}^{n}$ is called a domain of positivity
(D.P.), if there exists a non-singular symmetric matrix $S$, called a
characteristic of $D,$ such that $D$ is a maximal set so that $%
y^{t}Sx>0$ for all $y\in\overline{D}, y\neq 0$, and all $x\in D.$
\end{defi}
Clearly self-dual cones are domains of positivity with $S=I$, (the identity
matrix).

\begin{theo} Every D.P. is a cone.
\end{theo}

\begin{proof} We verify that all the conditions which define a cone are satisfied.
\begin{enumerate}
\item[(1)] If $x\in D, \lambda>0$, then $y^{t}S(\lambda x)=\lambda(y^{t}Sx)>0$ for
all $y\in\overline{D}, y\neq 0.$
\item[(2)] If $x, z\in D$, then $y^{t}S(x+z)=y^{\mathrm{t}}Sx+y^{t}Sz>0$ for all $%
y\in\overline{D}, y\neq 0.$
\item[(3)] Assume $D$ contains a line, i.e., $h+\lambda x\in D$, for all real $%
\lambda$. Taking $\lambda=0$ gives $h\in D$ so that $y^{t}S(h+\lambda x)>0$
for all $\lambda$. Thus, $y^{t}Sx=0$ for all $y\in D. \;D$ is open so that $%
x=0 $ and we do not have a line.
\end{enumerate}
\end{proof}

\begin{lem} Suppose $V$ is a cone and that $A$ is a regular matrix. Then $AV$
is a cone and $A^{t}$ maps $(AV)^{*}$ onto $V^{*}$
\end{lem}

\begin{proof} It is clear that $AV$ is a cone. Let $x\in(AV)^{*}$. Since $A%
\overline{V}=\overline{AV}$, then $\forall y\in\overline{V},$

$y\neq 0, \;A^{t}x\cdot y=x \cdot Ay>0$. Hence $A^{t}x\in V^{*}.$

By the first part of the proof since $(A^{t})^{-1}=(A^{-1})^{t}$ : $%
V^{*}\rightarrow(AV)^{*}, A^{t}$ is onto. 
\end{proof}

\begin{theo} Suppose $V$ is a D.P. If $A$ is regular, then $AV$ is also a
D.P.
\end{theo}
\begin{proof} 
\begin{equation*}
((A^{-1})^{t}SA^{-1})^{t}=(A^{-1})^{t}S^{t}((A^{-1})^{t})^{t}=(A^{-1})^{t}SA^{-1}
\end{equation*}
shows that $(A^{-1})^{t}SA^{-1}$ is symmetric. It is clearly regular.
Finally, 
\begin{equation*}
A V \stackrel{A^{-1}}{\longrightarrow} V \stackrel{S}{\longrightarrow} V^{*} \stackrel{\left(A^{-1}\right)^{t}}{\longrightarrow}(A V)^{*}\end{equation*}
shows
$$y^{t}[(A^{-1})^{t}SA^{-1}]x>0\; \text{ for all} \; x\in AV, y\in\overline{AV}.$$
\end{proof}

\section{$*$-function}

\begin{defi} (see \cite{Koeche}). For every $x\in V$, define $x^{*}=-($grad $%
\log\phi_{V}(x))^{t}.$
\end{defi}

\begin{theo} If $x\in V$, then $x^{*}\in V^{*}.$
\end{theo}
\begin{proof}

For every $y\in\overline{V}, y\neq 0,$

\begin{eqnarray*}
x^{*} \cdot y &=&-\left(\operatorname{grad} \log \phi_{V}(x)\right)^{t} \cdot y \\ 
&=&-\frac{1}{\phi_{V}(x)}\left(\operatorname{grad} \phi_{V}(x)\right)^{t} \cdot y \\ 
&=&-\frac{1}{\phi_{V}(x)}\left(\operatorname{grad} \int_{V^{*}} e^{-x \cdot z} d z\right)^{t} \cdot y \\ 
&=&\frac{1}{\phi_{V}(x)} \int_{V^{*}} e^{-x \cdot z}(z \cdot y) d z .
\end{eqnarray*}
Since $z\cdot y>0$, it follows that the last integral above is positive and so $%
x^{*} \cdot  y>0.$ Hence, $x^{*}\in V^{*}.$ 
\end{proof}

\begin{lem} If $V$ is a cone and $A$ is a regular map, then $\phi_{AV}(Ax)=|A|^{-1}\phi_{V}(x)$ .
\end{lem}
\begin{proof} 
\begin{eqnarray*}
\phi_{AV}(Ax)&=&\int_{(AV)^{*}}e^{-Ax\cdot y}dy\\
&=&\int_{(AV)^{*}} e^{-x\cdot A^{\mathrm{t}}y}dy\\
&=&\int_{V^{*}} e^{-x\cdot z}|A|^{-1}dz\\
&=&|A|^{-1}\phi_{V}(x).
\end{eqnarray*}
\end{proof}
In particular, taking $A=\lambda I$ with $\lambda>0$, then $AV=V$ and $
\phi_{V}=(\lambda x)=\lambda^{-n}\phi v(x)$ , i.e., $\phi_{V}$ is a
homogeneous function of degree $-n.$

\begin{lem} If $V$ is a cone and $A$ is a regular map, then $(Ax)^*=(A^{-1})^{t}x^{*}.$
\end{lem}
\begin{proof} Applying Lemma 2.11 we have,
\begin{eqnarray*}
(\textit{Ax})^* &=&-(\mathrm{g}\mathrm{r}\mathrm{a}\mathrm{d}%
_{Ax}\log\phi_{Av}(Ax))^t \\
&=&-(\mathrm{g}\mathrm{r}\mathrm{a}\mathrm{d}_{Ax}\log|A|^{-1}\phi_{V}(x))^{t}\\
&=&-(\mathrm{g}\mathrm{r}\mathrm{a}\mathrm{d}_{Ax}\log\phi_{V}(x))^{t}\\
&=&-[\mathrm{g}\mathrm{r}\mathrm{a}\mathrm{d}_{x}\log\phi_{V}(x)A^{-1}]^{t}\\
&=&(A^{-1})^{t}x^{*}.
\end{eqnarray*}
\end{proof}

\begin{defi} Let $V$ be a cone, define $K_{V}(x)=-\displaystyle \frac{\partial x^*}{\partial x}.$
\end{defi}
\begin{lem} If $V$ is a cone and $A$ is a regular map, then 
$$
K(Ax)=(A^{t})^{-1}K_V(x)A^{-1}.$$
\end{lem}
\begin{proof} Let $y=Ax.$ 
\begin{eqnarray*}
K_{AV}(Ax)=K_{AV}(y)&=&\frac{-\partial y^{*}}{\partial y}=\frac{
-\partial(Ax)^{*}}{\partial(Ax)}=\frac{-\left(A^{-1}\right)^{t} \partial x^{*}}{\partial x} \frac{\partial x}{\partial(A x)}\\
&=&- (A^{-1})^{t} \frac{\partial x^{*}}{\partial x} A^{-1}
=(A^{-1})^{t}K_{V}(x)A^{-1}.
\end{eqnarray*}
\end{proof}

\begin{lem} The Jacobian $K_{V}(x)=- \frac{\partial x^{*}}{\partial x}$ is non-singular and symmetric.
\end{lem}

\begin{proof}
\begin{eqnarray*}
x^{*} &=&-\left(\operatorname{grad} \log \phi_{V}(x)\right)^{t}=-\frac{1}{\phi_{V}(x)}\left(\operatorname{grad} \phi_{V}(x)\right)^{t} \\ 
&=& -\frac{1}{\phi_{V}(x)}\left[\begin{array}{c}\frac{\partial \phi_{V}(x)}{\partial x_{1}} \\ \vdots \\ \frac{\partial \phi_{V}(x)}{\partial x_{n}}\end{array}\right].
\end{eqnarray*}
So,
\begin{eqnarray*}
\frac{\partial x^{*}}{\partial x} &=&\left[\begin{array}{c}
-\frac{1}{\phi} \frac{\partial^{2} \phi}{\partial x_{1}^{2}}+\frac{1}{\phi^{2}}\left(\frac{\partial \phi}{\partial x_{1}}\right)^{2}, \ldots,-\frac{1}{\phi} \frac{\partial^{2} \phi}{\partial x_{n} \partial x_{1}}+\frac{1}{\phi^{2}}\left(\frac{\partial \phi}{\partial x_{n}}\right)\left(\frac{\partial \phi}{\partial x_{1}}\right) \\
\vdots \\
-\frac{1}{\phi} \frac{\partial^{2} \phi}{\partial x_{1} \partial x_{n}}+\frac{1}{\phi^{2}}\left(\frac{\partial \phi}{\partial x_{1}}\right)\left(\frac{\partial \phi}{\partial x_{n}}\right), \ldots,-\frac{1}{\phi} \frac{\partial^{2} \phi}{\partial x_{n}^{2}}+\frac{1}{\phi^{2}}\left(\frac{\partial \phi}{\partial x_{n}}\right)^{2}
\end{array}\right] \\
&=&-\frac{1}{\phi}\left[\begin{array}{c}
\frac{\partial^{2} \phi}{\partial x_{1}^{2}}, \ldots, \frac{\partial^{2} \phi}{\partial x_{n} \partial x_{1}} \\
\vdots \\
\frac{\partial^{2} \phi}{\partial x_{1} \partial x_{n}}, \ldots, \frac{\partial^{2} \phi}{\partial x_{n}^{2}}
\end{array}\right]+\frac{1}{\phi^{2}}\left[\begin{array}{c}
\left(\frac{\partial \phi}{\partial x_{1}}\right)^{2}, \ldots,\left(\frac{\partial \phi}{\partial x_{n}}\right)\left(\frac{\partial \phi}{\partial x_{1}}\right) \\
\vdots \\
\left(\frac{\partial \phi}{\partial x_{1}}\right)\left(\frac{\partial \phi}{\partial x_{n}}\right), \ldots,\left(\frac{\partial \phi}{\partial x_{n}}\right)^{2}
\end{array}\right].
\end{eqnarray*}
Note that $\frac{\partial x^*}{\partial x}$ is symmetric. Let $a\in {\mathbb R}^{n}, a\neq 0.$ 
\begin{eqnarray*}
a^{t} \frac{\partial x^{*}}{\partial x} a&=&-\frac{1}{\phi}\left[a_{1}, \ldots, a_{n}\right]\left[\begin{array}{c}\frac{\partial^{2} \phi}{\partial x_{1}^{2}}, \ldots, \frac{\partial^{2} \phi}{\partial x_{n} \partial x_{1}} \\ \vdots \\ \frac{\partial^{2} \phi}{\partial x_{1} \partial x_{n}}, \ldots, \frac{\partial^{2} \phi}{\partial x_{n}^{2}}\end{array}\right] \quad\left[\begin{array}{c}a_{1} \\ \vdots \\ a_{n}\end{array}\right]\\
&+&\frac{1}{\phi^{2}}\left[a_{1}, \ldots, a_{n}\right]\left[\begin{array}{c}\left(\frac{\partial \phi}{\partial x_{1}}\right)^{2}, \ldots,\left(\frac{\partial \phi}{\partial x_{n}}\right)\left(\frac{\partial \phi}{\partial x_{1}}\right) \\ \vdots \\ \left(\frac{\partial \phi}{\partial x_{1}}\right)\left(\frac{\partial \phi}{\partial x_{n}}\right), \ldots,\left(\frac{\partial \phi}{\partial x_{n}}\right)^{2}\end{array}\right]\left[\begin{array}{c}a_{1} \\ \vdots \\ a_{n}\end{array}\right]\\
&\equiv&-\frac{1}{\phi} \sum_{i, j} \frac{\partial^{2} \phi}{\partial x_{i} \partial x_{j}} a_{i} a_{j}+\frac{1}{\phi^{2}} \sum_{i, j}\left(\frac{\partial \phi}{\partial x_{i}}\right)\left(\frac{\partial \phi}{\partial x_{j}}\right) a_{i} a_{j}.
\end{eqnarray*}
Since 
\begin{equation*}
\frac{\partial^{2}\phi}{\partial x_{j}\partial x_{i}}=\int_{V^{*}}e^{-x\cdot
y}(y_{i}y_{j})dy,
\end{equation*}
\begin{eqnarray*} 
\sum_{i, j} \frac{\partial^{2} \phi}{\partial x_{i} \partial x_{j}} a_{i} a_{j} &=&\int_{V^{*}} e^{-x \cdot y}\left(\sum_{i, j} y_{i} y_{j} a_{i} a_{j}\right) d y \\ 
&=&\int_{V^{*}} e^{-x \cdot y}(a \cdot y)^{2} d y. 
\end{eqnarray*}
Also,
\begin{eqnarray*}
 \sum_{i, j}\left(\frac{\partial \phi}{\partial x_{i}}\right)\left(\frac{\partial \phi}{\partial x_{j}}\right) a_{i} a_{j} &=&\sum_{i, j}\left(\int_{V^{*}} e^{-x \cdot y}\left(-y_{i}\right) d y\right)\left(\int_{V^{*}} e^{-x \cdot y}\left(-y_{j}\right) d y\right) a_{i} a_{j} \\ &=&\sum_{i, j}\left(\int_{V^{*}} e^{-x \cdot y} a_{j} y_{j} d y\right)\left(\int_{V^{*}} e^{-x \cdot y} a_{i} y_{i} d y\right) \\ 
 &=&\sum_{j}\left(\sum_{i} \int_{V^{*}} e^{-x \cdot y}\left(a_{i} y_{i}\right) d y\right) \int_{V^{*}} e^{-x \cdot y} a_{j} y_{j} d y \\ &=&\left(\int_{V^{*}} e^{-x \cdot y}(a \cdot y) d y\right)^{2} .
\end{eqnarray*}

Since as a function of $y, e^{-\frac{x\cdot y}{2}}$ is not a constant multiple
of $e^{-\frac{x\cdot y}{2}} (a \cdot y)$ , then by Schwarz inequality, 
\begin{eqnarray*}
 \int_{V^{*}} e^{-x \cdot y}(a \cdot y) d y &<&\left(\int_{V^{*}} e^{-x \cdot y} d y\right)^{\frac{1}{2}}\left(\int_{V^{*}} e^{-x \cdot y}(a \cdot y)^{2} d y\right)^{\frac{1}{2}} \\ 
 &=&\left[\phi_{V}(x)\right]^{1 / 2}\left(\int_{V^{*}} e^{-x \cdot y}(a \cdot y)^{2} d y\right)^{1 / 2},
\end{eqnarray*}
so that 
\begin{equation*}
\frac{1}{\phi^{2}}\left(\int_{V^*}\ e^{-x\cdot y}(a\cdot y)dy\right)^{2}-\frac{1}{\phi}
\int_{V^*}\ e^{-x\cdot y}(a\cdot y)^{2}dy<0.
\end{equation*}
Hence,
$$a^{t}\displaystyle \frac{\partial x^{*}}{\partial x}a<0\; \text{for every} \;a\neq
0,$$
implying 
\begin{equation*}
|\frac{\partial x^{*}}{\partial x}|\neq 0.
\end{equation*}
\end{proof}

\begin{theo} (see \cite{Rothaus}). If $V$ is a cone, then $x^{*}=-($grad $
\log\phi_{V}(x))^{t}$ maps $V$ onto $V^{*}$ and satisfies the following
conditions:
\begin{enumerate}
\item[(i)] $\phi_{V}(x)\phi_{V^*}(x^{*})=c$
\item[(ii)] $|\displaystyle \frac{\partial x^*}{\partial x}|=c\phi_{V}^{2}(x)=c%
\phi_{V^{*}}^{-2}(x^{*})$
\item[(iii)] $(x^{*})^{*}=x$
\end{enumerate}
\end{theo}

\begin{proof} To show (i), take any $x, y\in V$. Since $V$ is homogeneous, there
exists $A\in G(V)$ such that $y=Ax$. So, by Lemma 2.11, 
\begin{equation*}
\phi_{V}(y)=\phi_{V}(Ax)=|A|^{-1}\phi_{V}(x)\ .
\end{equation*}
Also, 
\begin{equation*}
\phi_{V^{*}}(y^{*})=\phi_{V^{*}}((Ax)^{*})=\phi_{V^*}
((A^{-1})^{t}x^{*})=|(A^{-1})^t|^{-1}
\phi_{V^{*}}(x^{*})=|A|\phi_{V^*}(x^{*})\ .
\end{equation*}
Thus,
$$\phi_{V}(y)\phi_{V^*}(y^{*})=\phi_{V^*}(x)\phi_{V^{*}}(x^{*})$$
for every $x, y\in V,$
and 
\begin{equation*}
\phi_{V}(x)\phi_{V^*}(x^{*})=\mathrm{c}.
\end{equation*}
To prove (ii), by Lemma 2.13 we clearly have that if $x\in V$, then for
every $A\in G(V),$ 
\begin{equation*}
K_{V}(Ax)=(A^{-1})^{t}K_{V}(x)A^{-1}.
\end{equation*}
Taking determinants we get, 
\begin{equation*}
|K_{V}(Ax)|=|A|^{-2}|K_{V}(x)|.
\end{equation*}
Fix $x_{0}\in V$. Since $V$ is homogeneous, for any $y\in V,$ there exists $
A_{y}\in G(V)$ such that $A_{y}:x_{0}\rightarrow y$. So, 
\begin{equation*}
|K_{V}(y)|=|A_{y}|^{-2}|K_{V}(x_{0})|.
\end{equation*}
Also, by Lemma 2.11, $\phi_{V}(y)=\phi_{V}(A_{y}x_{0})=|A_{y}|^{-1}%
\phi_{V}(x_{0})$ . It follows that 
\begin{equation*}
|K_{V}(y)|=c\phi_{V}^{2}(y)
\end{equation*}
where 
\begin{equation*}
c=\frac{|K_{V}(x_{0})|}{\phi_{V}^{2}(x_{0})}
\end{equation*}
Now, using part (i),
$$|K_{V}(y)|=c\phi_{V}^{2}(y)=c\phi_{V^{*}}^{-2}(y^{*}),\; \text{for any}\;y\in V.$$

To prove (iii): Since from (i): 
$$\phi_{V}(x)\phi_{V^{*}}(x^{*})=c,$$
we have 
$$\begin{aligned} \log \phi_{V}(x) &+\log \phi_{V^{*}}\left(x^{*}\right)=\log c \\ & \Rightarrow \operatorname{grad}_{x^{*}} \log _{V}(x)+\operatorname{grad}_{x^{*}} \log \phi_{V^{*}}\left(x^{*}\right)=0 \\ & \Rightarrow\left(\frac{\partial x^{*}}{\partial x}\right)\left(\operatorname{grad}_{x^{*}} \log \phi_{V}(x)\right)^{t}+\left(\frac{\partial x^{*}}{\partial x}\right)\left(-x^{* *}\right)=0, \text { where } x^{* *}=\left(x^{*}\right)^{*} \\ & \Rightarrow\left(\operatorname{grad}_{x} \log \phi_{V}(x)\right)^{t}+\left(\frac{\partial x^{*}}{\partial x}\right)\left(-x^{* *}\right)=0.\end{aligned}$$
So,
\begin{equation}
x^{*}=(\displaystyle \frac{-\partial x^{*}}{\partial x})x^{**}.
\end{equation}

Since $\phi_{V}(x)$ is homogeneous of degree $-n$, we have from Euler's
formula, 
$$\displaystyle \sum_{i=1}^{n}x_{i}\frac{\partial\phi
}{\partial x_i}=-n\phi.$$
Therefore, 
\begin{eqnarray*}
\left(x^{*}\right)^{t} x &=&\left[\begin{array}{c}-\frac{\partial}{\partial x_{1}} \log \phi_{V}(x) \\ \vdots \\ -\frac{\partial}{\partial x_{n}} \log \phi_{V}(x)\end{array}\right]^{t}\left[\begin{array}{c}x_{1} \\ \vdots \\ x_{n}\end{array}\right] \\ 
&=&-x_{1} \frac{\partial}{\partial x_{1}} \log \phi_{V}(x)-\ldots-x_{n} \frac{\partial}{\partial x_{n}} \log \phi_{V}(x) \\ 
&=&-\frac{1}{\phi}\left(x_{1} \frac{\partial \phi}{\partial x_{1}}+\ldots+x_{n} \frac{\partial \phi}{\partial x_{n}}\right) \\ 
&=&-\frac{1}{\phi}(-n \phi)=n .
\end{eqnarray*}
Now, 
\begin{equation*}
n=(x^{*})^{t}x=\sum_{k=1}^{n}x_{k}^{*}x_{k}
\end{equation*}
implies 
$$\begin{aligned} 0 &=\frac{\partial}{\partial x_{i}} \sum_{k=1}^{n} x_{k}^{*} x_{k} \\ &=\sum_{k=1}^{n}\left[\left(\frac{\partial}{\partial x_{i}} x_{k}^{*}\right) x_{k}+x_{k}^{*} \frac{\partial x_{k}}{\partial x_{i}}\right] \\ &=\sum_{k=1}^{n}\left(\frac{\partial}{\partial x_{i}} x_{k}^{*}\right) x_{k}+x_{i}^{*}.\end{aligned}$$
Hence,
$$x_{i}^{*}=-\displaystyle \sum_{k=1}^{n}(\frac{\partial}{\partial x_{\dot{l}}} x_{k}^{*})x_{k}\; \text{for} \;i=1,2, \cdots n.$$

That is to say,
\begin{equation}
x^{*}=(-\displaystyle \frac{\partial x^{*}}{\partial x})x.
\end{equation}

Since the Jacobian $(\displaystyle \frac{\partial x^{*}}{\partial x})$ was
proved to be non-singular, by combining (2.1) and (2.2) we get that 
\begin{equation*}
(x^{*})^{*}=x,
\end{equation*}
and the mapping $x\rightarrow x^{*}$ is onto $V^{*}.$
\end{proof}

\begin{cor} $V$ is homogeneous if and only if $V^{*}$ is homogeneous.
\end{cor}

\begin{proof} Suppose $V$ is homogeneous. If $x,y\in V^{*}$, then $x^{*},
y^{*}\in V$ and so there exists $A\in G(V)$ such that $Ax^{*}=y^{*}.$ So, $
x^{*}=A^{-1}y^{*}= (A^ty)^*$ implying that $x=A^{t}y$ for $%
A^{t}\in G(V^{*})$. Hence, $V^{*}$ is homogeneous. The converse follows
from $V^{**}=V. $
\end{proof} 
The exposition to the end of the chapter is a
reorganization of results in \cite{Rothaus}.

\begin{lem} Assume that $V$ is a cone. If for some $x_{0}\in V, K_V(x_{0})\in
G(V\rightarrow V^{*}),$ where $G(V\rightarrow V^{*})$ is the group of
linear transformations mapping $V$ onto $V^{*}$, then for any $x\in V,
K_{V}(x)\in G(V\rightarrow V^{*}).$ 
\end{lem}

\begin{proof} Let $x\in V$. Let $A\in G(V)$ so that $Ax_{0}=x$. Let $v\in V$. Then 
\begin{equation*}
K_{V}(x)v=K_{V}(Ax_{0})v=(A^{-1})^{t}K_V(x_{0})A^{-1}v\in V^{*}.
\end{equation*}
Now to show that $K_{V}(x)$ is onto, let $y\in V^{*}$ It is obvious that $%
AK_{V}^{-1}(x_{0})A^{t}y\in V$ and that 
$$
K_{V}(x)[AK_{V}^{-1}(x_{0})A^{t}y]=(A^{-1})^{t}K_V(x_{0})A^{-1}[AK_{V}^{-1}(x_{0})A^{t}y]=y.$$
\end{proof}

\begin{lem} Let $V$ be a cone. Assume that $x_{0}, x_{0}^{*}$ and $0$ are on
a hne $L$. Then for any $x\in L$, we have $x^{*}\in L$. Further there exists 
$x_{1}\in L$ such that $x_{1}^{*}=x_{1}.$
\end{lem}

\begin{proof} For any $\lambda>0,$

\begin{equation}
K_{V}(\displaystyle \lambda x)=K_{V}(\lambda
Ix)=(\lambda^{-1}I)^{t}K_{V}(x)(\lambda I)^{-1}=\frac{1}{\lambda^{2}}%
K_{V}(x).
\end{equation}

Let $\lambda_{0}>0$, so that $x_{0}^{*}=\lambda_{0}x_{0}$. Then 
\begin{equation*}
(\lambda x_{0})^{*}=K_{V}(\lambda x_{0})(\lambda x_{0})=\frac{1}{\lambda^{2}}%
\lambda K_{V}(x_{0})x_{0}
\end{equation*}
\begin{equation*}
=\frac{1}{\lambda}x_{0}^{*}=\frac{\lambda_{0}}{\lambda}x_{0}=\frac{%
\lambda_{0}}{\lambda^{2}}(\lambda x_{0})\ .
\end{equation*}
This shows that for any $\lambda>0, (\lambda x_{0})^{*}$ is also on $L$.
Choosing $\lambda=\sqrt{\lambda_{0}}$, we see that $(\lambda
x_{0})^{*}=\lambda x_{0}.$
\end{proof}

\begin{theo} Let $V$ be a cone. Then there exists $x\in V$ so that $x=x^{*}$.
\end{theo}

\begin{proof} Consider the function $f(x)=\log\phi_{V}(x)$ defined on $\{\Sigma
x_{i}^{2}=n\}\cap V.$

Clearly $f(x)\rightarrow\infty$ as $x\rightarrow\partial V$. So $f(x)$
attains its minimum value at a point

$x^{0}\in\{\Sigma x_{i}^{2}=n\}\cap V$. At that point, using Lagrange
multipliers, we must have: 
\begin{equation*}
\sum(x_{i}^{0})^{2}=n
\end{equation*}
and
$$
\left\{\begin{array}{l}
\frac{\partial \log \phi_{V}(x)}{\partial x_{1}}+2 \lambda x_{1}^{0}=0 \\
\vdots \\
\frac{\partial \log \phi_{V}(x)}{\partial x_{n}}+2 \lambda x_{n}^{0}=0
\end{array}\right.
$$
That is,

$$\left\{\begin{array}{l}\left(x^{0}\right)_{1}^{*}=2 \lambda x_{1}^{0} \\ \vdots \\ \left(x^{0}\right)_{n}^{*}=2 \lambda x_{n}^{0}.\end{array}\right.$$

So, 
\begin{equation*}
\sum_{\dot{l}}(x^{0})_{i}^{*}\cdot x_{i}^{0}=2\lambda n.
\end{equation*}
Since $(x^{0})^{*}\cdot x^{0}=n, \displaystyle \lambda=\frac{1}{2}$ Thus, 
\begin{equation*}
(x^{0})^{*}=x^{0}.
\end{equation*}
\end{proof}

\begin{theo} Let $V$ be a D.P. with characteristic $S$. Recall that $S$ is
a symmetric matrix. Then for every $x\in V,$ 
\begin{equation*}
K_{V}(S^{-1}x^{*})=SK_{V}^{-1}(x)S.
\end{equation*}
\end{theo}
\begin{proof} Let $v=S^{-1}x^{*}\in V$. Then $x^{*}=Sv$ and $%
x=(Sv)^{*}=(S^{-1})v^{*}.$ Thus, 
\begin{equation*}
K_{V}(S^{-1}x^{*})=K_{V}(v)=\frac{-\partial v^{*}}{\partial v}=\frac{%
\partial v^{*}}{\partial x}(-\frac{\partial x}{\partial x^{*}})\frac{%
\partial x^{*}}{\partial v}=SK_{V}^{-1}(x)S.
\end{equation*}
\end{proof}

\begin{theo} Let $V$ be a domain of positivity with positive definite
characteristic $S.$ Then there exists $p\in V$ such that $Sp=p^{*}$
\end{theo}
\begin{proof} We will look for the infimum of $\log\phi_{V}(x)$ subject to $
x^{t}Sx=n$. Since $S$ is positive definite (and symmetric), there exists a
unique positive definite (and symmetric) $P$ such that $P^{2}=S$. Thus, 
$$\begin{aligned} x^{t} S x=n & \Leftrightarrow x^{t} P^{t} P x=n \\ & \Leftrightarrow(P x)^{t} P x=n \\ & \Leftrightarrow\|P x\|=\sqrt{n},\end{aligned}$$
implying that such $x$'s form a bounded set, say $B.$

Consider $\log\phi_{V}(x)$ over $B\cap V$. Since $\phi_{V}(x)\rightarrow%
\infty$ as $x\rightarrow\partial V, \log\phi v(x)$ attains a minimum at an
interior point of $V$ subject to $x^{\mathrm{t}}Sx=n.$
Call this point $p$. At $p$ the following equations must hold: 
\begin{equation*}
\frac{\partial}{\partial x_{k}}\{\log\phi_{V}(x)+\lambda(x^{t}Sx-n)\}=0,\
k=1,2,\cdots, n,
\end{equation*}
for a suitable value of $\lambda.$ 
$$\begin{aligned} \frac{\partial}{\partial x_{k}}\left(x^{t} S x\right) &=\frac{\partial}{\partial x_{k}}\left(\sum_{i=1}^{n} \sum_{j=1}^{n} S_{i j} x_{i} x_{j}\right) \\ &=\sum_{i=1}^{n} S_{i k} x_{i}+\sum_{j=1}^{n} S_{k j} x_{j}=2 \sum_{j=1}^{n} S_{k j} x_{j} \\ &=2(S x)_{k}. \end{aligned}$$
We will have $p^{*}=2\lambda Sp$ and 
\begin{equation*}
n=(p^{*})^{\mathrm{t}}p=2\lambda(Sp)^{t}p=2\lambda p^{t}S^{t}p=2\lambda p^{%
\mathrm{t}}Sp=2\lambda n.
\end{equation*}
Therefore, $\displaystyle \lambda=\frac{1}{2}$ and $p^{*}=Sp. $
\end{proof}

\begin{theo} Let $V$ be a domain of positivity with positive definite
characteristic $S$. Let $Sp=p^{*}$. Then $K_{V}(p)=S$, and for any $x\in V,
K_{V}(x)$ is a positive definite characteristic.
\end{theo}

\begin{proof} From Theorem 2.19, for any $x\in V,$ 
\begin{equation*}
S^{-1}K_{V}(x)S^{-1}K_{V}(S^{-1}x^{*})=I.
\end{equation*}
So, in particular, since $Sp=p^{*},
[S^{-1}K_{V}(p)]^{2}=I$. Since $S$ and $Kv(p)$ are positive definite, $%
K_{V}(p)=S.$

Given $x\in V$, by Lemma 2.15 $K_{V}$ is a characteristic. Further, given
any $a\in {\mathbb R}^{n}, a\neq 0$, let $Ap=x, A\in G(V)$; then

$$a^{t}K_{V}(x)a=a^{\mathrm{t}
}(A^{-1})^{t}K_{V}(p)A^{-1}a=(A^{-1}a)^{t}K_{V}(p)(A^{-1}a)=(A^{-1}a)^{t}S(A^{-1}a)>0. 
$$
\end{proof}

\begin{theo} Let $V$ be a domain of positivity with characteristic $S$.
Then for any $a, b\in V$, there exists $x\in V$ so that $K(x)a=Sb.$
\end{theo}

\begin{proof} Consider the minimum of $x^{*}\cdot a$ for $x\in V$ subject to $%
x^{t}Sb=1.$

From Lemma 2.8, $x^{t}Sb\geq r(Sb)$ for all $x\in\Sigma_{\overline{V}} $ where $r(Sb)=$ distance of $Sb$ from $\partial V^{*}$, we know
that for all $x\in\overline{V}, x^{t}Sb\geq r(Sb)|x|$. Thus,

$$|x|\displaystyle \leq\frac{x^{t}Sb}{r(Sb)}=\frac{1}{r(Sb)}\; \text{for all} \; 
x\in\{x\in\overline{V}:x^{t}Sb=1\}.$$

We claim that the function $x^{*}$ . $a\rightarrow\infty$ as $%
x\rightarrow\partial V, x\in V$. Let us show first that the minimum of $%
\phi_{V}(y)$ subject to $x^{*}\cdot y=n$ is $\phi_{V}(x)$ .

From $x^{*}\cdot y\geq r(x^{*})$ for all $y\in\Sigma_{\overline{V}}$, we
know for all $y\in V, x^{*}\cdot y\geq r(x^{*})|y|$, so that $|y|
\displaystyle \leq\frac{x^{*}\cdot y}{r(x^{*})}=\frac{n}{r(x^{*})}$ and
therefore, the set $\{y\in V\ :\ x^{*}\cdot y=n\}$ is bounded.

Since $\log\phi_{V}(y)\rightarrow\infty$ as $y\rightarrow\partial V$, there
must be a point $y_{0}\in V$, where $\log\phi_{V}$ has its minimum. We have: 
\begin{equation*}
\frac{\partial}{\partial y_{j}}\{\log\phi_{V}(y)-\lambda(x^{*}\cdot
y-n)\}=0,\; \;j=1,\ \cdots n
\end{equation*}
holds at $y=y_{0}$. It follows that $y_{0}^{*}=\lambda x^{*}.$ But then, from 
$n=y_{0}^{*}. y_{0}=\lambda (x^{*}\ .\ y_{0})=\lambda n, \lambda=1$. So, we
get $\phi_{V}(y)\geq\phi_{V}(x)$ for all $y$ so that $x^{*}\cdot y=n.$

Given $a\in V$, let $y=\displaystyle \frac{n}{a\cdot x^{*}}z\in V$. Then $%
x^{*}\cdot y=n$, and thus $\displaystyle \phi_{V}(\frac{n}{a\cdot x^{*}}%
z)\geq\phi_{V}(x).$ This implies that $(\displaystyle \frac{n}{a\cdot x^{*}}%
)^{-n}\phi_{V}(z)\geq\phi_{V}(x),$ i.e.,

$$\left(a \cdot x^{*}\right)^{n} \geq \frac{\phi_{V}(x)}{\phi_{V}(z)} n^{n} \quad \text{for all} \;\; x, z \in V.$$

Since $\phi_{V}(x)\rightarrow\infty$ as $x\rightarrow\partial V$, we know
that $x^{*}\cdot a\rightarrow\infty$ as $x\rightarrow\partial V$ where $a\in
V.$

Finally, there exists a point $x_{0}\in V$ where $x^{*}\cdot a$ is minimal
subject to $x^{\mathrm{t}}Sb=1$. At this point

$$\frac{\partial}{\partial x_{i}}\left\{x^{*} \cdot a-\lambda\left(x^{t} S b-1\right)\right\}=0 \quad j=1,2, \ldots, n\;  \text{holds at} \; x=x_{0}.$$

It follows that $K_{V}(x_{0})a=\lambda Sb$. Since $K_{V}(x_{0})a, Sb$ are
both points of $V^{*}, \lambda$ is positive. Therefore, by (2.3), $K_{V}(
\displaystyle \sqrt{\lambda}x_{0})a=\frac{1}{\lambda}K_{V}(x_{0})a=Sb.$
\end{proof}

\begin{cor} Let $V$ be a domain of positivity with positive definite
characteristic. Then for any $a\in V, d\in V^{*},$ there exists $x\in V$
such that $K(x)a=d.$
\end{cor}

\begin{proof} If $d\in V^{*}, S^{-1}d\in V$, and so there exists $x\in V$ so that 
$$ K(x)a=S(S^{-1}d)=d.$$
\end{proof}

\begin{theo} Let $V$ be a domain of positivity with positive definite
characteristic. Let $a, b\in V$. Then for some constant $c>0,$ 
\begin{equation*}
c\phi_{V}(a+b)=\phi_{V}(a)\phi_{V}(b)\phi_{V^{*}}(a^{*}+b^{*})\ .
\end{equation*}
\end{theo}

\begin{proof} By Corollary 2.5, there exists $x\in V$ so that $K(x)a=b^{*}$. Then, $
b=(K(x)a)^{*}= (K_{V}^{-1}(x))a^{*}$ or $K_{V}(x)b=a^{*}.$

Thus, $K_{V}(x)(a+b)=a^{*}+b^{*}$ so that 

$$\phi_{V^{*}}(a^{*}+b^{*})=\phi_{V^{*}}(K(x)(a+b))=|K_{V}(x)|^{-1}\phi_{V}(a+b).$$
 Similarly, from $K(x)a=b^{*}$, we get 
 $$
|Kv(x)|^{-1}\phi_{V}(a)=\phi_{V^{*}}(b^{*})= \displaystyle \frac{c}{\phi_{V}(b)}$$
so that $|K_{V}(x)|=\displaystyle \frac{\phi_V(a) \phi_V(b)}{c}$ and $c\phi_{V}(a+b)=\phi_{V}(a)\phi_{V}(b)\phi_{V^*}(a^{*}+b^{*})$ . 
\end{proof}

\begin{theo} Let $V$ be a domain of positivity with positive definite
characteristic. Let $a,b\in V$. Then 
\begin{equation*}
(a+b)^{*}=a^{*}-K(a)(a^{*}+b^{*})^{*}.
\end{equation*}
\end{theo}
\begin{proof} From Theorem 2.23, we have 
\begin{equation*}
\log
c+\log\phi_{V}(a+b)=\log\phi_{V}(a)+\log\phi_{V}(b)+
\log_{V^{*}}(a^{*}+b^{*}).
\end{equation*}
Taking the grad with respect to the variable $a$, we get 
\begin{equation*}
\mathrm{g}\mathrm{r}\mathrm{a}\mathrm{d}_{a}\log\phi_{V}(a+b)=\mathrm{g}
\mathrm{r}\mathrm{a}\mathrm{d}_{a}\log\phi_{V}(a)+\mathrm{g}\mathrm{r}
\mathrm{a}\mathrm{d}_{a}*\log\phi_{V^{*}}(a^{*}+b^{*})(\frac{\partial
a^{*}}{\partial a})\ .
\end{equation*}
Thus, 
\begin{equation*}
(a+b)^{*}=a^{*}-K(a)(a^{*}+b^{*})^{*}.
\end{equation*}
\end{proof}

Given a cone $V$ we define a partial ordering in ${\mathbb R}^{n}$.

\begin{defi} $x \underset{V}{<} y$ if and only if $y-x\in V.$
\end{defi}

\begin{theo} If $V$ is a domain of positivity with positive definite
characteristic, then
$$x \underset{V}{<} Z \quad \text{iff} \quad z^{*}  \underset{V^*}{<} x^{*}.$$
\end{theo}

\begin{proof} Let $y=z-x\in V$. From Theorem 2.24 we have that 
\begin{equation*}
(x+y)^{*}=x^{*}-K(x)(x^{*}+y^{*})^{*}.
\end{equation*}
Hence, 
\begin{equation*}
x^{*}-z^{*}=K(x)(x^{*}+(z-x)^{*})^{*}\in V^{*},
\end{equation*}
i. e., 
\begin{equation*}
 z^{*} \underset{V^*}{<} x^{*}.
\end{equation*}
Conversely, if $ z^{*} \underset{V^*}{<} x^{*}$, then $(x^{*})^{*}\underset{(V^*)^*}{<} (z^{*})^{*}$,
i.e., $x \underset{V}{<} z. $
\end{proof}

\chapter{FUNCTIONS AND OPERATORS ON CONES}

We begin with a review of properties of $V$-homogeneous functions and
integral operators on cones and then prove the boundedness of such operators
on some weighted $L^{p}$ spaces.

\section{Functions Defined on Cones}

\begin{defi} Let $V$ be a cone and let $a, b\in V$. The cone interval $
\langle a, b\rangle$ is defined: 
\begin{equation*}
\langle a, b \rangle=\{x \in V:a\underset{V}{<}x\underset{V}{<}b\}.
\end{equation*}
\end{defi}
In ${\mathbb R}^{2}, \langle a, b\rangle$ is a parallelogram having vertices at $a$ and $b$
and edges parallel to the boundaries of the cone.

\begin{defi} Let $V$ be a cone. The function $\triangle_{V}(x)$ is
defined: 
\begin{equation*}
\triangle_{V}(x)=\int_{\langle 0,x\rangle}dy,\ x\in V,
\end{equation*}
i.e., $\triangle_{V}(x)$ is the measure of the cone interval $\langle 0,
x\rangle.$
\end{defi}

\begin{defi} A function $f$ : $V\rightarrow {\mathbb R}^{+}$ is said to be $V$%
-homogeneous of order $\delta$ if $f(Ax)=|A|^{\delta}f(x)\; \forall A \in
G(V)$ where $|A|$ denotes the absolute value of the determinant of the
matrix $A.$
\end{defi}

It was shown earlier that $\phi_{V}(Ax)=|A|^{-1}\phi_{V}(x)$ for any $A\in
G(V)$ . Thus, $\phi_{V}(x)$ is $V$-homogeneous of order $-1$. Also $%
\triangle_{V}(x)$ is $V$-homogeneous of order 1 since for any $A\in G(V)$ , 
\begin{equation*}
\triangle_{V}(Ax)=\int_{\langle
0,Ax\rangle}dy=\int_{\langle0,x\rangle}|A|dz=|A|\triangle_{V}(x)\ .
\end{equation*}
Note that if $f$ is $V$-homogeneous of any order, then $f$ is either
identically $0$ on $V$ or $f(x)\neq 0$ for all $x\in V.$

 It is also known that $V$-homogeneous functions are also
homogeneous in the usual sense, since if $f$ is $V$-homogeneous of order $%
\delta$, then for $A=\lambda I$ with $\lambda>0,$ 
\begin{equation*}
f(\lambda x)=f(Ax)=|A|^{\delta}f(x)=\lambda^{n\delta}f(x)\ ,
\end{equation*}
i.e., $f$ is homogeneous of degree $n\delta.$

\begin{theo}(see \cite{Ostrogorski}). If $V$ is a homogeneous cone, then all $V$-homogeneous
functions are, up to a multiplicative constant, powers of $\triangle_{V}(x)$.
\end{theo}

\begin{proof} Assume first that $f(x)$ is $V$-homogeneous of order $0$. For any $x,
y\in V$, there exists $A\in G(V)$ so that $Ax=y$. Hence, 
\begin{equation*}
f(y)=f(Ax)=f(x)\ .
\end{equation*}
So, $f(x)=c.$

Now if $f(x)$ is $V$-homogeneous of order $\delta$, then 
$F(x)= \frac{f(x)}{\triangle_{V}^{\delta}(x)}$ is $V$
-homogeneous of order $0$, so that

$F(x)=c$ and $f(x)=c\triangle_{V}^{\delta}(x)$.
\end{proof} 

From the above theorem, 
\begin{equation*}
\phi_{V}(x)=c\triangle_{V}^{-1}(x)=\frac{c}{\triangle_{V}(x)}.
\end{equation*}
We can now translate the properties of $\phi_{V}(x)$ into properties of $
\triangle_{V}(x)$ .
\begin{enumerate}
\item[(i)] \begin{equation} \triangle_{V}(x)\triangle_{V^{*}}(x^{*})=c\end{equation}
\item[(ii)] \begin{equation}\label{5} |K(x)|=|\displaystyle \frac{\partial x^{*}}{\partial x}|=c\Delta_{V}^{-2}(x)=c
\triangle_{V^{*}}^{2}(x^{*})\end{equation}
\item[(iii)] \begin{equation} \triangle_{V}(x) \;\text{ is continuous and} \; \triangle_V(x)>0\:  \text{for every} \; x\in V\end{equation}
\item[(iv)] \begin{equation}\triangle_{V}(x)\rightarrow 0
\;\text{as}\; x\rightarrow\partial V\end{equation}
\end{enumerate}

\section{Integral Operators on Cones}

In this section, we will be considering integral operators of the form

\begin{equation}
Kf(x)=\displaystyle \int_{V}k(x,\ y)f(y)dy,\ x\in V
\end{equation}
where $k:V\times V\rightarrow {\mathbb R}^{+}$ and $f:V\rightarrow {\mathbb R}^{+}.$

\begin{defi} The kernel $k(x,y)$ : $V\times V\rightarrow {\mathbb R}^{+}$ is said
to be $V\times V$-homogeneous of order $\beta$ if

$$k(Ax,\ Ay)=|A|^{\beta}k(x,\ y)$$ for all $A\in G(V).$
\end{defi}
\begin{defi} For a positive function $f$
defined on $V$, Hardy's operator is defined as 
\begin{equation*}
Hf(x)=\int_{\langle 0,x\rangle}f(y)dy.
\end{equation*}
Also, Laplace's operator is defined as 
\begin{equation*}
Lf(x)=\int_{V}e^{-x^{*}\cdot y}f(y)dy.
\end{equation*}
\end{defi}
Since $x\underset{V}{<}y$ if and only if $Ax\underset{V}{<} Ay$ for all $A\in G(V),$ it can
easily be seen that the kernel of Hardy's operator

$$k(x,\ y)=\left\{
\begin{array}{ll}
1 & \text{if} \quad y\underset{V}{<}x \\ 
0 &\text{for other} \;  y\in V
\end{array}
\right.$$
is $V\times V$-homogeneous of order $0$. Also, the kernel of Laplace's
operator $k(x,\ y)=e^{-x^{*}\cdot y}$ is $V\times V$-homogeneous of order $0$, since 
\begin{equation*}
k(Ax,\ Ay)=e^{-(Ax)^{*}\cdot Ay}=e^{-(A^{-1})^{t}x^*\cdot Ay}=e^{-x^{*}\cdot y}=k(x,\ y)\ .
\end{equation*}
In one dimension, Hardy's operator has been shown to be $L^{p}$-continuous
on weighted $L^{p}$ spaces with weights $x^{\gamma}$: (see for example \cite{Zygmund}).

If $1\leq p<\infty$ and $\gamma<p-1$, then
\begin{equation}
\displaystyle \int_{0}^{\infty}\left(\int_{0}^{x}f(y)dy\right)^{p}x^{\gamma-p}dx\leq
c\int_{0}^{\infty}f^{p}(x)x^{\gamma}dx,
\end{equation}
where $c$ is a constant independent of $f.$

In subsequent works, (see for example \cite{Muckenhoupt}) the weights $x^{\gamma}$ were
replaced by a larger class of functions and necessary and sufficient
conditions for Hardy's operator to become $L^{p}$-continuous were found. If $
1<p<\infty$, then
\begin{equation}
\left( \int_{0}^{\infty}|u(x)\int_{0}^{x}f(y)dy|^{p}dx\right)^{1/p}\leq c\left(\int_{0}^{\infty}|f(x)v(x)|^{p} dx \right)^{1/p}
\end{equation}
 if and only if
\begin{equation}
\displaystyle \sup_{r>0}\left(\int_{r}^{\infty}|u(x)|^{p}dx\right)^{1/p}\left(\int_{0}^{r}|v(x)|^{-p^{\prime }}dx\right)^{1/p'}=A<\infty,
\end{equation}
and $A\leq c\leq Ap^{1/p}p^{\prime 1/p^{\prime }}.$ There also has been a $%
[L^{p},\ L^{q}]$ generalization of the above by J. S. Bradley: (see \cite{Bradley}).

If $u, v$ are non-negative and $1\leq p\leq q\leq\infty$, then
\begin{equation}
\displaystyle\left(\int_{0}^{\infty}(u(x)\int_{0}^{x}f(y)dy)^{q}dx\right)^{1/q}\leq
c\left(\int_{0}^{\infty}(f(x)v(x))^{p}dx\right)^{1/p}
\end{equation}
if and only if
\begin{equation}
\displaystyle \sup_{r>0}\left(\int_{r}^{\infty}u^{\mathrm{q}}(x)dx\right)^{1/q}\left(
\int_{0}^{r}v^{-p^{\prime }}(x)dx\right)^{1/p^{\prime }}=A<\infty .
\end{equation}

Furthermore, $A\leq c\leq Ap^{1/p}p^{\prime1/p^{\prime }}$ for $1<p\leq q<\infty$
and $A=c$ if $\mathrm{p}=1$ or $q=\infty.$

The generalization to ${\mathbb R}^{n}$ of Hardy's inequality for self-dual cones is
due to T. Ostrogorski \cite{Ostrogorski}.

\begin{theo} (Ostrogorski). Let $V$ be a homogeneous, self-dual cone in $%
{\mathbb R}^{n}$ and

$1\leq p<\infty$. Assume that $k$ is a $V\times V$-homogeneous kernel of
order $0$ and $k(x,\ y)=k(y^{*},\ x^{*})$ for all $x, y\in V$. If the
integral $K\triangle_{V}^{\alpha}(x)$ is convergent for some $\alpha\in R$,
then
\begin{equation}\label{14}
\displaystyle \int_{V}(K(f(x))^{p}\triangle_{V}^{\gamma-p}(x)dx\leq
c\int_{V}f^{\mathrm{p}}(x)\triangle_{V}^{\gamma}(x)dx,
\end{equation}
where $\gamma=-\alpha p-1.$
\end{theo}

As special cases of $K$ on self-dual cones, Hardy's and Laplace's operators
were considered in \cite{Ostrogorski}. It was shown that if

$\displaystyle \sigma_{0}(V)=\inf\{\alpha\in R\ :\
\int_{\Sigma_{V}}\triangle_{V}^{\alpha}(t)dt<\infty\}$ and $\displaystyle %
\sigma(V)=\max\{-1,\ \sigma_{0}\},$

then for $\alpha>\sigma(V)$ , both $H\triangle_{V}^{\alpha}(x)$ and $%
L\triangle_{V}^{\alpha}(x)$ are finite, so that:

\begin{cor}  If $\gamma<-\sigma(V)p-1$, then

\begin{equation}
\displaystyle \int_{V}(Hf(x))^{p}\triangle_{V}^{\gamma-p}(x)dx\leq
c\int_{V}f^{p}(x)\triangle_{V}^{\gamma}(x)dx
\end{equation}
and

\begin{equation}\displaystyle \int_{V}(Lf(x))^{p}\triangle_{V}^{\gamma-p}(x)dx\leq
c\int_{V}f^{p}(x)\triangle_{V}^{\gamma}(x)dx.
\end{equation}
\end{cor} 

\section{Norm Inequalities on Cones}

The condition in Ostrogorski's theorem that the underlying cone be self-dual
is too restrictive. We extend her results to more general cones.

\begin{theo} Let $V$ be a cone in ${\mathbb R}^{n}$ and $1\leq p\leq q<\infty$. Assume
that the kernel $k(x,\ y)$ : $V\times V\rightarrow {\mathbb R}^{+}$ of operator $K$ is 
$V\times V$-homogeneous of order $\beta$. If for some $\delta,\gamma\in R,$

\begin{equation}\label{17} 
K\displaystyle \triangle_{V}^{\delta}(x)=\int_{V}k(x,\
y)\triangle_{V}^{\delta}(y)dy<\infty\end{equation}
and

\begin{equation}\label{18}
\displaystyle \int_{V}k^{q/p}(x,\
y)\triangle_{V}^{\delta-q+(\delta+\beta+1)q/p^{l}}(x)dx<\infty,\end{equation}
then

\begin{equation}\label{19}
\displaystyle \left(\int_{V}\triangle_{V}^{\gamma-q}(x)(Kf(x))^{q}dx\right)^{1/q}\leq
c\left(\int_{V}f^{p}(x)\Delta_{V}^{\beta p+(\gamma+1)p/q-1}(x)dx\right)^{1/p}.
\end{equation}
\end{theo}
\begin{proof} Using H\"{o}lder's inequality, we have 
$$
\begin{aligned}
& \int_{V} \Delta_{V}^{\gamma-q}(x)(K f(x))^{q} d x \\
=& \int_{V} \Delta_{V}^{\gamma-q}(x)\left(\int_{V} k^{1 / p}(x, y) f(y) \Delta_{V}^{-\delta / p^{\prime}}(y) k^{1 / p^{\prime}}(x, y) \Delta_{V}^{\delta / p^{\prime}}(y) d y\right)^{q} d x \\
\leq & \int_{V} \Delta_{V}^{\gamma-q}(x)\left(\int_{V} k(x, y) f^{p}(y) \Delta_{V}^{-\delta(p-1)}(y) d y\right)^{q / p}\left(\int_{V} k(x, y) \Delta_{V}^{\delta}(y) d y\right)^{q / p^{\prime}} d x \\
\end{aligned}
$$
where $\displaystyle \frac{1}{p}+\frac{1}{p}=1.$

Now, $K\triangle_{V}^{\delta}(x)$ is finite and for every $A\in G(V),$ 
$$\begin{aligned}
K \Delta_{V}^{\delta}(A x) &=\int_{V} k(A x, y) \Delta_{V}^{\delta}(y) d y \\
&=\int_{V} k(A x, A u) \Delta_{V}^{\delta}(A u)|A| d u \\
&=|A|^{\beta+\delta+1} \int_{V} k(x, u) \Delta_{V}^{\delta}(u) d u
\end{aligned}$$
so that $K\triangle_{V}^{\delta}(x)$ is $V$-homogeneous of order $%
\delta+\beta+1$. Hence, 
\begin{equation*}
\int_{V}k(x,\ y)\triangle_{V}^{\delta}(y)dy=c\triangle_{V}^{\delta+\beta+1}(x)\ .
\end{equation*}

We therefore have: 
\begin{eqnarray*}
 &&\int_{V} \Delta_{V}^{\gamma-q}(x)(K f(x))^{q} d x \\
&& \leq c \int_{V} \Delta_{V}^{\gamma-q+(\delta+\beta+1) q / p^{\prime}}(x)\left(\int_{V} k(x, y) f^{p}(y) \Delta_{V}^{-\delta(p-1)}(y) d y\right)^{q / p} d x.
\end{eqnarray*}
Since $q/p\geq 1$, we can apply Minkowski's integral inequality to the last
integral. Thus, 
\begin{eqnarray*}
&&\int_{V}\triangle_{V}^{\gamma-q}(x)(Kf(x))^{q}dx\\
&& \leq c\left(\int_{V}f^{p}(y)\triangle_{V}^{-\delta(p-1)}(y)\left(\int_{V}k^{q/p}(x,y)\triangle_{V}^{\gamma-q+(\delta+\beta+1)q/p^{\prime
}}(x)dx\right)^{p/q}dy\right)^{q/p}.
\end{eqnarray*}
Again, 
\begin{equation*}
\int_{V}k^{q/\mathrm{p}}(x,\
y)\triangle_{V}^{\gamma-q+(\delta+\beta+1)q/p^{\prime }}(x)dx<\infty
\end{equation*}
for each $y$, and can readily be shown to be $V$-homogeneous of order $
\gamma+\beta q+\delta q/p^{\prime }-q/p+1.$ It therefore equals $
c\triangle_{V}^{\gamma+\beta q+\delta q/p^{\prime }-q/p+1}(y)$ .
Substituting then, 
\begin{equation*}
\int_{V}\triangle_{V}^{\gamma-q}(x)(Kf(x))^{q}dx\leq c\left(\int_{V}f^{p
}(y)\triangle_{V}^{(\gamma+1)p/q+\beta p-1}(y)dy\right)^{q/p}.
\end{equation*}
Raising both sides to $(1/q)$-th power gives the result. 
\end{proof}

Note
that if $V=V^{*}, p=q$ and $\beta=0$, we obtain (\ref{14}).\\

The following theorem corresponds to the case $q=\infty$ in Theorem 3.3.

\begin{theo} Let $V$ be a cone and $1\leq p<\infty$. Assume that the kernel $%
k(x,\ y)$ : $V\times V\rightarrow {\mathbb R}^{+}$ is $V\times V$-homogeneous of order 
$\beta$. Assume also that

\begin{equation}
K\displaystyle \triangle_{V}^{\delta}(x)=\int_{V}k(x,\
y)\triangle_{V}^{\delta}(y)dy<\infty, \; \text{for some}\; \delta,
\end{equation}
and
\begin{equation}
\displaystyle \mathrm{e}\mathrm{s}\mathrm{s}\sup_{x\in V}k(x,y)\triangle_{V}^{(\delta+\beta)(p-1)-1}(x)<\infty.
\end{equation}
Then

\begin{equation}
\displaystyle   \mathrm{e}\mathrm{s}\mathrm{s}\sup_{x\in
V}\triangle_{V}^{-1}(x)(Kf(x))\leq c\left(\int_{V}f^{p}(x)\triangle_{V}^{\beta p-1}(x)dx\right)^{1/p}.
\end{equation}
\end{theo}
\begin{proof} Let $g(y)=\displaystyle \mathrm{e}\mathrm{s}\mathrm{s}\sup_{x\in
V}k(x,\ y)\triangle_{V}^{\langle\delta+\beta)(p-1)-1}(x).$ For every $A\in
G(V)$,
\begin{eqnarray*}
g(Ay)&=& \mathrm{e}\mathrm{s}\mathrm{s}\sup_{x\in V}k(x, Ay)\Delta_{V}^{(\delta+\beta)(p-1)-1}(x)\\
&=&\displaystyle  \mathrm{e}\mathrm{s}\mathrm{s}\sup_{x\in V}k(Ax,Ay)\triangle_{V}^{(\delta+\beta)(p-1)-1} (Ax) \\
&=&|A|^{\delta(p-1)+\beta p-1}\mathrm{e}\mathrm{s}\mathrm{s}\sup_{x\in V}k(x,y)\triangle_{V}^{(\delta+\beta)(p-1)-1}(x).
\end{eqnarray*}

So, $g(y)$ is $V$-homogeneous of order $\delta(p-1)+\beta p-1$ and
therefore, 
\begin{equation*}
g(y)=c\triangle_{V}^{\delta(p-1)+\beta p-1}(y)\ .
\end{equation*}
Now, by H\"{o}lder's inequality

\begin{eqnarray*}
\displaystyle \triangle_{V}^{-1}(x)(Kf(x))&=&\Delta_{V}^{-1}(x) \int_{V} k^{1 / p}(x, y) \Delta_{V}^{-\delta / p^{\prime}}(y) f(y) k^{1 / p^{\prime}}(x, y) \Delta_{V}^{\delta / p^{\prime}}(y) d y\\
&\leq&\triangle_{V}^{-1}(x)\left(\int_{V}k(x,\
y)\triangle_{V}^{-\delta(p-1)}(y)f^{p}(y)dy\right)^{1/p}\left(\int_{V}k(x,\
y)\triangle_{V}^{\delta}(y)dy\right)^{1/p^{\prime }}\\
&=&c\triangle_{V}^{(\delta+\beta+1)1/p^{\prime }-1}(x)\left(\int_{V}k(x,\
y)\triangle_{V}^{-\delta(p-1)}(y)f^{p}(y)dy\right)^{1/p}\\
&=&c\left(\int_{V}k(x,y)\triangle_{V}^{-p+(\delta+\beta+1)(p-1)}(x)\triangle_{V}^{-\delta(p-1)}(y)f^{p}(y)dy\right)^{1/p}\\
&\leq& c\left(\int_{V}(\triangle_{V}^{-\delta(p-1)+\beta p-1}(y)\triangle_{
V}^{\delta(p-1)}(y)f^{p}(y)dy\right)^{1/p}\\
&=&c\left(\int_{V}\triangle_{V}^{\beta p-1}(y)f^{p}(y)dy\right)^{1/p}.
\end{eqnarray*}
\end{proof}

The conclusion of Theorem 3.4 can also be proved under somewhat different
conditions.

\begin{theo} Let $V$ be a cone and $1\leq p<\infty$. Assume that the kernel $
k$ is $V\times V$ homogeneous of order $\beta$. Assume also that

\begin{equation}\label{23}
g(x)=\left(\displaystyle \int_{V}k^{p^{\prime }}(x,\ y)\triangle_{V}^{(1-\beta)
p^{\prime}-1}(y)dy\right)^{1/p^{\prime}}<\infty.
\end{equation}

Then

\begin{equation}\label{24}
\displaystyle \mathrm{e}\mathrm{s}\mathrm{s}\sup_{x\in
V}\triangle_{V}^{-1}(x)(Kf(x))\leq c\left(\int_{V}f^{p}(y)\triangle_{V}^{\beta
p-1}(y)dy\right)^{1/p}.
\end{equation}
\end{theo}
\begin{proof} Applying H\"{o}lder's inequality, we get

\begin{eqnarray*}
\displaystyle \triangle_{V}^{-1}(x)(Kf(x))&=&\triangle_{V}^{-1}(x)%
\int_{V}k(x,\ y)\Delta_{V}^{-\beta+1/p}(y)f(y)\triangle_{V}^{\beta-1/p}(y)dy\\
&\leq&\triangle_{V}^{-1}(x)\left(\int_{V}f^{p}(y)\triangle_{V}^{\beta
p-1}(y)dy\right)^{1/p}\left(\int_{V}k^{p^{\prime }}(x,
y)\triangle_{V}^{(1-\beta)p^{\prime }-1}(y)dy\right)^{1/p^{\prime }}.
\end{eqnarray*}
Since $g(x)<\infty$ and is $V$-homogeneous of order 1, we have $
g(x)=c\triangle_{V}(x)$ . Therefore, 
\begin{equation*}
\triangle_{V}^{-1}(x)(Kf(x))\leq c\left(\int_{V}f^{p}(y)\triangle_{V}^{\beta
p-1}(y)dy\right)^{1/p}.
\end{equation*}
\end{proof}
The following theorem corresponds to the case $p=\infty$ and $q=\infty$ in
Theorem 3.3.

\begin{theo} Let $V$ be a cone. Assume that the kernel
$k(x,\ y)$ : $V\times V\rightarrow {\mathbb R}^{+}$ is $V\times V$-homogeneous of
order $\beta$. Assume also that

\begin{equation}\label{25}
K\displaystyle \triangle_{V}^{\delta}(x)=\int_{V}k(x,\
y)\triangle_{V}^{\delta}(y)dy<\infty\;\; \text{for some} \;\delta.
\end{equation}

Then

\begin{equation}
\displaystyle \mathrm{e}\mathrm{s}\mathrm{s}\sup_{x\in
V}\triangle_{V}^{-1-\delta-\beta}(x)(Kf(x))\leq c\;\mathrm{e}\mathrm{s}\mathrm{%
s}\sup_{x\in V}f(x)\triangle_{V}^{-\delta}(x).
\end{equation}
\end{theo}
\begin{proof}
\begin{eqnarray*}
\triangle_{V}^{-1-\delta-\beta}(x)(Kf(x))&=&\triangle_{V}^{-1-\delta-\beta}(x)\int_{V}k(x,\ y)\triangle_{V}^{\delta}(y)f(y)\triangle_{V}^{-\delta}(y)dy\\
&\leq&\triangle_{y}^{-1-\delta-\beta}(x)\left(\int_{V}k(x,\
y)\triangle_{V}^{\delta}(y)dy\right)\mathrm{e}\mathrm{s}\mathrm{s}\sup_{y\in
V}f(y)\triangle_{V}^{-\delta}(y)\\
&=&c\;\mathrm{e}\mathrm{s}\mathrm{s}\sup_{y\in V}f(y)\triangle_{V}^{-\delta}(y).
\end{eqnarray*}
\end{proof}

\section{Norm Inequalities for Dual Operators}

Another reasonable generalization of Ostrogorski's theorem is to consider
kernels 
$$k(x,\ y):V^{*}\times V\rightarrow {\mathbb R}^{+}.$$

\begin{defi} The kernel $k(x,\ y)$ : $V^{*}\times V\rightarrow {\mathbb R}^{+}$ is
said to be $V^{*}\times V$-homogeneous of order $\beta$ if

$$k ((A^{t})^{-1}x, Ay)=|A|^{\beta}k(x,y)$$
 for all $A\in G(V)$ .
\end{defi}

The following shows that all results concerning such kernels can be deduced
from results on kernels $$k:V\times V\rightarrow {\mathbb R}^{+}.$$

\begin{theo} Let $V$ be a cone. If $K\triangle_{V^*}^{\delta}(x)$ is finite
and the kernel $k:V^{*}\times V\rightarrow {\mathbb R}^{+}$ is $V^{*}\times V$%
-homogeneous of order $\beta,$ then there exists a constant $c$
such that 
\begin{equation*}
\int_{V}k(x,\ y)\triangle_{V}^{\delta}(y)dy=c\int_{\langle 0,x^{*}\rangle}k(x,y)\triangle_{V}^{\delta}(y)dy.
\end{equation*}
\end{theo}

\begin{proof} Let $\psi(x)$ be the right-hand side integral. For every $A\in G(V)$, 
\begin{eqnarray*}
\psi((A^{t})^{-1}x)&=&\int_{\langle 0,((A^t)^{-1}x)^*\rangle}k((A^{t})^{-1}x,\ y)\triangle_{V}^{\delta}(y)dy\\
&=&\int_{\langle 0,Ax^*\rangle}k((A^{t})^{-1}x,Au)\triangle_{V}^{\delta}(Au)d(Au)\\
&=&\int_{\langle 0,x^{*}\rangle}|A|^{\beta}k(x,u)|A|^{\delta}\triangle_{V}^{\delta}(u)|A|du\\
&=&|A|^{\beta+\delta+1}\psi(x)\ .
\end{eqnarray*}
Hence, $\psi(x)$ is $V^{*}$-homogeneous of order $-(\beta+\delta+1)$ .
Similarly,
$$K\displaystyle \triangle_{V^{*}}^{\delta}(x)=\int_{V}k(x,\
y)\triangle_{V}^{\delta}(y)dy$$ is $V^{*}$-homogeneous of order $
-(\beta+\delta+1)$ . It therefore follows that 
\begin{equation*}
\int_{V}k(x,\ y)\triangle_{V}^{\delta}(y)dy=c\int_{\langle 0,x^*\rangle}k(x,\
y)\triangle_{V}^{\delta}(y)dy.
\end{equation*}
\end{proof}

\begin{defi} We define an operator $S$ mapping measurable functions on $V$
onto measurable functions on $V^{*}$ by: $Sf(x)=f(x^{*})$ where $x\in V^{*}.$
\end{defi}

\begin{theo} 
\begin{equation*}
\int_{V^{*}}(Sf(x))^{q}\triangle_{V^{*}}(x)dx=c\int_{V}f^{q}(y)\triangle_{V}^{-\delta-2}(y)dy.
\end{equation*}
\end{theo}

\begin{proof} Using (\ref{5}), we obtain, 
\begin{eqnarray*}
\int_{V^{*}}(Sf(x))^{q}\triangle_{V^{*}}^{\delta}(x)dx&=&c\int_{V}f^{q}(x^{*})
\triangle_{V}^{-\delta}(x^{*})|\frac{\partial x}{\partial x^{*}}|dx^{*}\\
&=&c\int_{V}f^{q}(x^{*})\triangle_{V}^{-\delta}(x^{*})
\triangle_{V}^{-2}(x^{*})dx^{*}\\
&=&c\int_{V}f^{q}(y)\triangle_{V}^{-\delta-2}(y)dy.
\end{eqnarray*}
\end{proof}
Applying Theorem 3.8, we obtain the following versions of Theorems 3.3-3.6
for $V^{*}\times V$ operators.

\begin{theo}  Let $V$ be a cone in ${\mathbb R}^{n}$ and $1\leq p\leq q<\infty$. Assume
that $k(x,\ y)$ is $V^{*}\times V$-homogeneous of order $\beta$. If

\begin{equation}
\displaystyle \int_{V}k(x, y)\triangle_{V}^{\delta}(y)dy<\infty
\end{equation}

and

\begin{equation}\label{28}
 \int_{V^{*}} k^{q/p}(x,y)\triangle_{V^{*}}^{-\gamma+q-2-(\delta+\beta+1)q/p^{\prime}}(x)dx<\infty,
\end{equation}
then
\begin{equation}
\left( \int_{V^*} (Kf(x))^{q}\triangle_{V^{*}}^{-\gamma+q-2}(x)dx\right)^{1/q}\leq
c\left(\int_{V}f^{p}(y)\triangle_{V}^{\beta p+(\gamma+1)p/q-1}(y)dy\right)^{1/p}.
\end{equation}
\end{theo}
\begin{proof} We will apply Theorem 3.3. Define $\tilde{k}(u,\ y)$ : $V\times
V\rightarrow {\mathbb R}^{+}$ by 
\begin{equation*}
\tilde{k}(u,\ y)=k(u^{*}, y).
\end{equation*}
Clearly, $\tilde{k}$ is $V\times V$-homogeneous of order $\beta$ if and only
if $k$ is $V^{*}\times V$-homogeneous of order $\beta.$

Then 
\begin{equation*}
SKf(u)=\int_{V}\tilde{k}(u,\ y)f(y)dy.
\end{equation*}
By hypothesis, 
\begin{equation*}
\int_{V}\tilde{k}(u,\ y)\triangle_{V}^{\delta}(y)dy=\int_{V}k(u^{*},\
y)\triangle_{V}^{\delta}(y)dy<\infty.
\end{equation*}
Also, by Theorem 3.8 and (\ref{28}),
\begin{eqnarray*}
\displaystyle \int_{V}\tilde{k}^{q/p}(u, y)\triangle_{V}^{\gamma-q+(\delta+\beta+1)q/p^{\prime }}(u)du
&=&c\int_{V^*}(S\tilde{k})^{q/p}(x,y)\triangle_{V^*}^{-\gamma+q-2-(\delta+\beta+1)q/p\prime}(x)dx\\
&=&c\int_{V^{*}}k^{q/p}(x,y)\triangle_{V^*}^{-\gamma+q-2-(\gamma+\beta+1)q/p^{\prime }}(x)dx<\infty.
\end{eqnarray*}
Note that $S(\displaystyle \int_{V}\tilde{k}(u,\ y)f(y)dy)=Kf(u)$ . So by
Theorems 3.3 and 3.8, 
\begin{eqnarray*}
\left(\int_{V^*}(Kf(u))^{q}\triangle_{V^{*}}^{-\gamma+q-2}(u)du\right)^{1/q}&=&c(\int_{V}(
\int_{V}\tilde{k}(u,\ y)f(y)dy)^{q}\triangle_{V}^{\gamma-q}(u)du)^{1/q}\\
&\leq& c(\int_{V}f^{p}(x)\triangle_{V}^{\beta p+(\gamma+1)p/q-1}(x)dx)^{1/p}.
\end{eqnarray*}
\end{proof}
The following theorem corresponds to the case $q=\infty$ in Theorem 3.9.

\begin{theo} Let $V$ be a cone and $1\leq p<\infty$. Assume that the kernel 
$k(x,\ y)$ : $V^{*}\times V\rightarrow {\mathbb R}^{+}$ is $V^{*}\times V$-homogeneous
of order $\beta$. If

\begin{equation}
\displaystyle \int_{V}k(x, y)\triangle_{V}^{\delta}(y)dy<\infty\;\quad \text{for some} \;
\delta,
\end{equation}
and
\begin{equation}
\displaystyle \mathrm{e}\mathrm{s}\mathrm{s}\sup_{x\in V^{*}}k(x,\
y)\triangle_{V^{*}}^{-(\delta+\beta)(p-1)+1}(x)<\infty,
\end{equation}

then
\begin{equation}
\displaystyle \mathrm{e}\mathrm{s}\mathrm{s}\sup_{x\in
V^{*}}\triangle_{V^{*}}(x)(Kf(x))\leq c(\int_{V}f^{\mathrm{p}}(y)\triangle_{V}^{\beta p-1}(y)dy)^{1/p}.
\end{equation}
\end{theo}
\begin{proof}  The proof is immediate if we note that $g(y)=\displaystyle \mathrm{e}\mathrm{s}\mathrm{s}\sup_{x\in V^*}k(x,\
y)\triangle_{V^{*}}^{-(\delta+\beta)(p-1)+1}(x)$ is $V$-homogeneous of order 
$\delta(p-1)+\beta p-1$, and follow the steps used in the proof of Theorem
3.4. 
\end{proof}
The conclusion of Theorem 3.10 also follows from different conditions.

\begin{theo} Let $V$ be a cone and $1\leq p<\infty$. Assume that the kernel 
$k$ is $V^{*}\times V$-homogeneous of order $\beta.$

If 
\begin{equation*}
\int_{V}k^{p^{\prime }}(x,\ y)\triangle_{V}^{(1-\beta\rangle p^{\prime}-1}(y)dy<\infty,
\end{equation*}
then 
\begin{equation*}
\mathrm{e}\mathrm{s}\mathrm{s}\sup_{x\in
V^{*}}\triangle_{V^{*}}(x)(Kf(x))\leq c(\int_{V}f^{p}(y)\triangle_{V}^{\beta
p-1}(y)dy)^{1/p}.
\end{equation*}
\end{theo}
\begin{proof} Let $g(x)=(\displaystyle \int_{V}k^{p^{\prime }}(x,y)\triangle_{V}^{(1-\beta)p^{\prime }-1}(y)dy)^{1/p^{\prime }}.$ It can
easily be seen that $g(x)$ is $V^{*}$-homogeneous of order $-1$ and so
equals $c\triangle_{V^{*}}^{-1}(x)$ . The proof is now immediate if we
follow the same line of argument as that in the proof of Theorem 3.5.
\end{proof}
\begin{theo} Let $V$ be a cone. Assume that the kernel $k(x,\ y)$ : $%
V^{*}\times V\rightarrow {\mathbb R}^{+}$ is $V^{*}\times V$-homogeneous of order $%
\beta.$
If

$$\displaystyle \int_{V}k(x,\ y)\triangle_{V}^{\delta}(y)dy<\infty\quad \text{for some}\;
\delta,$$
then 
\begin{equation*}
\mathrm{e}\mathrm{s}\mathrm{s}\sup_{x\in
V^*}\triangle_{V^{*}}^{1+\delta+\beta}(x)(Kf(x))\leq c\, \mathrm{e}\mathrm{s}%
\mathrm{s}\sup_{y\in V}f(y)\triangle_{V}^{-\delta}(y) .
\end{equation*}
\end{theo}

\begin{proof} Let $K\displaystyle \triangle_{V^{*}}^{\delta}(x)=\int_{V}k(x,y)\triangle_{V}^{\delta}(y)dy$. It was seen earlier that $%
K\triangle_{V^{*}}^{\delta}(x)$ is $V^{*}$-homogeneous of order $%
-(\beta+\delta+1)$ . The proof can now be readily obtained following the
same lines as those in the proof of Theorem 3.6.
\end{proof}

\section{ Applications}

In this section we will apply the results of our general theorems to some
special operators: Riemann-Liouville's and Weyl's fractional integral
operators and Laplace's operator. Riemann-Liouville's and Weyl's
inequalities will be proved on domains of positivity with positive definite
characteristic. Laplace's inequalities, however, will be proved in the more
general setting of homogeneous cones. Note that Hardy's inequalities are
special cases of Riemann-Liouville's. See \cite{Martin-Reyes} for some related recent
results on Riemann-Liouville's fractional integral operators in $R.$

\begin{defi} The Riemann-Liouville operator is defined: 
\begin{equation*}
R_{r}f(x)=\frac{1}{\Gamma(r)}\int_{\langle
0,x\rangle}\triangle_{V}^{r-1}(x-t)f(t)dt,\ x\in V.
\end{equation*}
The Weyl operator is defined: 
\begin{equation*}
W_{r}f(x)=\frac{1}{\Gamma(r)}\int_{\langle
x,\infty\rangle}\triangle_{V}^{r-1}(t-x)f(t)dt,\ x\in V,
\end{equation*}
where $r\geq 1.$
\end{defi}
Note that by Fubini's theorem
for any non-negative measurable functions $f$ and $g$ defined on $V,$

\begin{equation}\label{33} \int_{V}(R_{r}g)\cdot f=\int_{V}g\cdot(W_{r}f).
\end{equation}

Note also that the cases $r=1$ in Theorems 3.13 and 3.14 below provide
generalizations of the well-known Hardy's inequalities from $R$ to ${\mathbb R}^{n}$.

\begin{theo} (Riemann-Liouville's Inequalities). Let $V$ be a domain of
positivity with positive definite characteristic. If $1\leq p\leq q<\infty$
and $\gamma<-\sigma(V) q / p^{\prime}-\sigma\left(V^{*}\right)+q(1 / p- r+1)-2$, then

\begin{equation}\label{34}
\left(\displaystyle \int_{V}\triangle_{V}^{\gamma-q}(x)(R_{r}f(x))^{q}dx\right)^{1/q}\leq
c\left(\int_{V}f^{p}(x)\triangle_{V}^{(r-1)p+(\gamma+1)p/q-1}(x)dx\right)^{1/p
}.
\end{equation}

If $1\leq p<\infty$, and $\alpha<2-(1+\sigma(V))/p^{\prime }-r$, then

\begin{equation}\label{35}
\displaystyle \mathrm{e}\mathrm{s}\mathrm{s}\sup_{x\in
V}\left(\triangle_{V}^{-1+\alpha}(x)R_{r}(f\triangle_{V}^{-\alpha})(x)\right)\leq
c\; \left(\int_{V}f^{p}(x)\triangle_{V}^{(r-1)p-1}(x)dx\right)^{1/p}.
\end{equation}

If $\alpha<1-r-\sigma(V),$ then

\begin{equation}\label{36}
\underset{x \in V}{\operatorname{ess} \sup }\left(\Delta_{V}^{-1+\alpha}(x) R_{r}\left(f \Delta_{V}^{-\alpha}\right)(x)\right) \leq c\underset{x \in V}{\operatorname{ess}}\left(f(x) \Delta_{V}^{r-1}(x)\right).
\end{equation}
\end{theo}

\begin{proof} We prove (\ref{34}) as an application of Theorem 3.3. The kernel 
\begin{equation*}
k(x,t)=\frac{1}{\Gamma(r)}\triangle_{V}^{r-1}(x-t)\chi_{\langle
0,x\rangle}(t)
\end{equation*}
is $V\times V$-homogeneous of order $r-1.$

Now, since $x \underset{V}{<} y$ implies $\triangle_{V}(x)<\triangle_{V}(y)$ for all $x,
y\in V$, we have, 
\begin{eqnarray*}
R_{r}\triangle_{V}^{\delta}(x)&=&\frac{1}{\Gamma(r)}
\int_{\langle 0,x)}\triangle_{V}^{r-1}(x-t)\triangle_{V}^{\delta}(t)dt\\
&\leq&\frac{1}{\Gamma(r)}\int_{\langle
0,x\rangle}\triangle_{V}^{r-1}(x)\triangle_{V}^{\delta}(t)dt\\
&=&\frac{\triangle_{\mathrm{t}}^{r-1}(x)}{\Gamma(r)}\int_{\langle
0,x\rangle}\triangle_{V}^{\delta}(\mathrm{t})dt
\end{eqnarray*}
and the last integral is finite if $\delta>\sigma(V)$ (see \cite{Ostrogorski}). Moreover,
since $V$ is a domain of positivity with positive definite characteristic, $t \underset{V}{<} x$ if and only if $x^* \underset{V^*}{<} t^8$. Therefore,

\begin{eqnarray*}
\Gamma^{q / p}(r) \int_{V} k^{q / p}(x, t) \Delta_{V}^{\gamma-q+(\delta+r) q / p^{\prime}}(x) d x &=&\int_{\langle t, \infty\rangle} \Delta_{V}^{(r-1) q / p}(x-t) \Delta_{V}^{\gamma-q+(\delta+r) q / p^{\prime}}(x) d x \\
& \leq &\int_{\langle t, \infty\rangle} \Delta_{V}^{(r-1) q / p+\gamma-q+(\delta+r) q / p^{\prime}}(x) d x \\
&=&c \int_{\left\langle 0, t^{*}\right\rangle} \Delta_{V^{*}}^{-(r-1) q / p-\gamma+q-(\delta+r) q / p^{\prime}-2}\left(x^{*}\right) d x^{*}.
\end{eqnarray*}

Again, the integral on the right side is finite if 
\begin{equation*}
-(r-1)q/p-\gamma+q-(\delta+r)q/p^{\prime}-2 > \sigma(V^*).
\end{equation*}
Now, the condition on $\gamma$ in hypothesis ensures the existence of $\delta\in \mathbb{R}$ such that $-(r-1)q/p-\gamma+q-(\delta+r)q/p^{\prime}-2 > \sigma(V^*).$ Applying Theorem 3.3, the proof is complete.

To prove (\ref{35}), we apply Theorem 3.5 to the kernel

\begin{equation}\label{37}
k(x,\displaystyle \ t)=\frac{1}{\Gamma(r)}\chi_{\langle 0,x\rangle}(t)
\triangle_{V}^{\alpha}(x)\triangle_{V}^{-\alpha}(t)%
\triangle_{V}^{r-1}(x-t)
\end{equation}

Clearly, $k(x,\ t)$ is $V\times V$-homogeneous of order $r-1$. Also, 
\begin{eqnarray*}
&&\int_{\langle 0,x\rangle}\triangle_{V}^{\alpha p^{\prime}}(x)
\triangle_{V}^{-\alpha p^{\prime}}(t)\triangle_{V}^{(r-1)p^{\prime
}}(x-t)\triangle_{V}^{(2- r)p^{\prime}-1}(t)dt\\
&& \hspace{1cm} \leq\triangle_{V}^{\alpha p^{\prime }+(r-1)p^{\prime }}(x)\int_{\langle
0,x\rangle}\triangle_{V}^{(2-r-\alpha)p^{\prime }-1}(t)dt,
\end{eqnarray*}
and the last integral is finite since $\alpha<2-(1+\sigma(V))/p^{\prime }-r$.
The conditions of Theorem 3.5 are satisfied and we get (35).

To prove (\ref{36}), we apply Theorem 3.6 to the kernel $k(x,t)$ of (\ref{37}). Let $%
\delta=1-r$ in (25). We have, 
\begin{equation*}
\frac{1}{\Gamma(r)}\int_{\langle
0,x\rangle}\triangle_{V}^{\alpha}(x)\triangle_{V}^{-\alpha}(t)%
\triangle_{V}^{r-1}(x-t)\triangle_{V}^{1-r}(t)dt\leq\frac{%
\triangle_{V}^{\alpha+r-1}(x)}{\Gamma(r)}\int_{\langle
0,x\rangle}\triangle_{V}^{1-r-\alpha}(t)dt.
\end{equation*}
Now, since $\alpha<1-r-\sigma(V)$, the last integral is finite, and the
proof is complete. \end{proof}

\begin{theo} (Weyl's Inequalities). Let $V$ be a domain of positivity with
positive definite characteristic.

If $1\leq p\leq q<\infty$ and $\gamma>\sigma(V)+\sigma(V^{*})q/p^{\prime
}+2q-q/p$, then

\begin{equation}\label{38}
\left(\displaystyle \int_{V}\triangle_{V}^{\gamma-q}(x)(W_{r}f(x)\right)^{q}dx)^{1/q}%
\leq
c\left(\int_{V}f^{p}(x)\triangle_{V}^{(r-1)p+(\gamma+1)p/q-1}(x)dx\right)^{1/p}.
\end{equation}

If $1\leq p<\infty$ and $\alpha>\sigma(V^{*})/p^{\prime }+2-1/p$, then

\begin{equation}\label{39}
\displaystyle \mathrm{e}\mathrm{s}\mathrm{s}\sup_{x\in
V}\left(\triangle_{V}^{-1+\alpha}(x)W_{r}(f\triangle_{V}^{-\alpha})(x)\right)\leq
c\left(\int_{V}f^{p}(x)\triangle_{V}^{(r-1)p-1}(x)dx\right)^{1/p}.
\end{equation}

If $\alpha>2+\sigma(V^{*})$ , then

\begin{equation}\label{40}
\displaystyle \mathrm{e}\mathrm{s}\mathrm{s}\sup_{x\in
V}(\triangle_{V}^{-1+\alpha}(x)W_{f}(f\triangle_{V}^{-\alpha})(x))\leq c%
\mathrm{e}\mathrm{s}\mathrm{s}\sup_{x\in V}(f(x)\triangle_{V}^{r-1}(x)).
\end{equation}
\end{theo}
\begin{proof} We will use the duality of the Weyl operator and the
Riemann-Liouville operator, (\ref{33}), to prove the theorem.

Let $g(x)$ be any non-negative function defined on $V$ with $\displaystyle 
\int_{V}g^{q'}(x)dx=1$, where $\displaystyle \frac{1}{q'}+\frac{1}{
q}=1$. Applying H\"{o}lder's inequality, we have 
\begin{eqnarray*}
 \int_{V} \Delta_{V}^{\gamma / q-1}(x) g(x)\left(W_{r} f(x)\right) d x &=& \int_{V} f(x) R_{r}\left(g \Delta_{V}^{\gamma / q-1}\right)(x) d x \\
& &\hspace{-3.5cm}\leq\left(\int_{V} f^{p}(x) \Delta_{V}^{(r-1) p+(\gamma+1) p / q-1}(x) d x\right)^{1 / p} \times \\
&&\hspace{-2.5cm} \left(\int_{V} \Delta_{V}^{-(r-1) p^{\prime}-(\gamma+1) p^{\prime} / q+p^{\prime}-1}(x)\left(R_{r}\left(g \Delta_{V}^{\gamma / q-1}\right)(x)\right)^{p^{\prime}} d x\right)^{1 / p^{\prime}}.
\end{eqnarray*}

We now apply Theorem $3.13$ to the last integral. Note that $q \geq p$ iff $p^{\prime} \geq q^{\prime}$; also $\gamma>\sigma(V)+\sigma\left(V^{*}\right) q / p^{\prime}+2 q-q / p$ implies that
$$
\left(-(r-1) p^{\prime}-(\gamma+1) p^{\prime} / q+p^{\prime}-1\right)+p^{\prime}<-\sigma(V) p^{\prime} / q-\sigma\left(V^{*}\right)+p^{\prime}\left(1 / q^{\prime}-r+1\right)-2.
$$
 Therefore, using (\ref{34}), we see that
the last integral is mayorized by 
\begin{equation*}
c\left(\int_{V}\left(g(x) \Delta_{V}^{\gamma/q-1}(x)\right)^{q^{\prime}} \Delta_{V}^{2 q^{\prime}-(\gamma+1) q^{\prime} / q-1}(x) d x\right)^{1 / q^{\prime}}=c\left(\int_{V} g^{q^{\prime}}(x) d x\right)^{1 / q^{\prime}}=c.
\end{equation*}
Hence, 
\begin{equation*}
\int_{V}\triangle_{V}^{\gamma/q-1}(x)g(x)(W_{r}f(x))dx\leq
c(\int_{V}f^{P}(x)\triangle_{V}^{(r-1)_{p}+(\gamma+1)p/q-1}(x)dx)^{1/p}.
\end{equation*}
Taking supremum over all $g\geq 0$ so that $\displaystyle 
\int_{V}g^{q^{\prime }}(x)dx=1$, we get 
\begin{equation*}
\left(\int_{V}\triangle_{V}^{\gamma/q-1}(x)(W_{r}f(x))^{q}dx\right)^{1/q}\leq
c\left(\int_{V}f^{p}(x)\triangle_{V}^{(r-1)p+(\gamma+1)p/q-1}(x)dx\right)^{1/p}.
\end{equation*}
To prove (\ref{39}), let $g(x)$ be any non-negative function on $V$ so that $%
\displaystyle \int_{V}g(x)dx=1$. We have: 
\begin{eqnarray*}
\int_{V}\triangle_{V}^{-1+\alpha}(x)g(x)W_{r}(f\triangle_{V}^{-\alpha})(x)dx &=&\int_{V}f(x)\triangle_{V}^{-\alpha}(x)R_{r}(g\triangle_{V}^{-1+\alpha})(x)dx\\
&&\hspace{-3.5cm}\leq\left(\int_{V}f^{p}(x)\triangle_{V}^{(r-1)p-1}(x)dx\right)^{1/p}\times\\
&&\hspace{-2.5cm}\left(\int_{V} \Delta_{V}^{-(r-1) p^{\prime}+\left(p^{\prime}-1\right)-\alpha p^{\prime}}(x)\left(R_{r}\left(g \Delta_{V}^{-1+\alpha}\right)(x)\right)^{p^{\prime}} d x\right)^{1 / p^{\prime}} .
\end{eqnarray*}
 Now, the condition $\alpha>\sigma(V^{*})/p^{\prime }+2-1/p$ implies that 
\begin{equation*}
-(r-1)p^{\prime }+(p^{\prime }-1)-\alpha p^{\prime }+p^{\prime}<-\sigma(V^{*})+p^{\prime }(1-r+1)-2
\end{equation*}
and hence using (\ref{34}), the last integral is majorized by 
\begin{equation*}
 c \int_{V} g(x) \Delta_{V}^{-1+\alpha}(x) \Delta_{V}^{(r-1)+\left(-(r-1) p^{\prime}+\left(p^{\prime}-1\right)-\alpha p^{\prime}+p^{\prime}+1\right) / p^{\prime}-1}(x) d x 
= c \int_{V} g(x) d x=c.
\end{equation*}
Therefore, 
\begin{equation*}
\int_{V}\triangle_{V}^{-1+\alpha}(x)g(x)W_{r}(f\triangle_{V}^{-\alpha})(x)dx%
\leq c\left(\int_{V}f^{p}(x)\triangle_{V}^{(r-1)p-1}(x)dx\right)^{1/p}.
\end{equation*}
Taking supremum over all $g\geq 0$ so that $\displaystyle \int_{V}g(x)dx=1$
completes the proof of (\ref{39}).

Finally, to prove (\ref{40}), let $g(x)$ be such that $\displaystyle %
\int_{V}g(x)dx=1$. We have, 
\begin{eqnarray*}
 \int_{V} \Delta_{V}^{-1+\alpha}(x) g(x) W_{r}\left(f \Delta_{V}^{-\alpha}\right)(x) d x 
&=& \int_{V} f(x) \Delta_{V}^{-\alpha}(x) R_{r}\left(g \Delta_{V}^{-1+\alpha}\right)(x) d x \\
&&\hspace{-3.5cm} \leq\operatorname{ess} \sup _{x \in V}\left(f(x) \Delta_{V}^{r-1}(x)\right) \int_{V} \Delta_{V}^{-\alpha-r+1}(x) R_{r}\left(g \Delta_{V}^{-1+\alpha}\right)(x) d x.
\end{eqnarray*}
By Fubini's theorem, 
\begin{eqnarray*}
\int_{V}\triangle_{V}^{-\alpha-r+1}(x)R_{r}(g\triangle_{V}^{-1+\alpha})(x)dx\\
&&\hspace{-3.5cm}=\frac{1}{\Gamma(r)}\int_{V}g(t)\triangle_{V}^{-1+\alpha}(t)\left(\int_{\langle t ,\infty
\rangle}\triangle_{V}^{r-1}(x-t)\triangle_{V}^{-\alpha-r+1}(x)dx\right)dt.
\end{eqnarray*}
Now, since $\alpha>2+\sigma(V^{*})$ , it can be shown that 
\begin{equation*}
\int_{\langle t,\infty\rangle}\triangle_{V}^{-\alpha-r+1}(x)\triangle_{V}^{r-1}(x-t)dx< \infty
\end{equation*}
and is $V$-homogeneous of order $-\alpha+1;$ it therefore equals $
c\triangle_{V}^{-\alpha+1}(t)$ . Hence, 
\begin{equation*}
\int_{V}\triangle_{V}^{-\alpha-r+1}(x)R_{r}(g\triangle_{V}^{-1+\alpha})(x)dx
\end{equation*}
\begin{equation*}
=c\int_{V}g(t)\triangle_{V}^{-1+\alpha-\alpha+1}(t)dt=c.
\end{equation*}
Therefore, 
\begin{equation*}
\int_{V}\triangle_{V}^{-1+\alpha}(x)g(x)W_{r}(f\triangle_{V}^{-\alpha})(x)dx%
\leq c\;\mathrm{e}\mathrm{s}\mathrm{s}\sup_{x\in
V}(f(x)\triangle_{V}^{r-1}(x))\ .
\end{equation*}
Taking supremum over all $g\geq 0$ with $\displaystyle \int_{V}g^{q^{\prime
}}(x)dx=1$, gives us (\ref{40}). \end{proof}

 Note that Theorem 3.14 could also have been proved by applying the
results of Theorems 3.3, 3.4 or 3.5 and 3.6.

\begin{theo}(Laplace's Inequalities). Let $V$ be a cone in ${\mathbb R}^{n}.$

If $1 \leq p \leq q<\infty$ and $\gamma<-\sigma(V) q / p^{\prime}-\sigma\left(V^{*}\right)+q / p-2$, then

\begin{equation}\label{41}
\left(\displaystyle \int_{V}\triangle_{V}^{\gamma-q}(x)(Lf(x))^{q}dx\right)^{1/q}\leq
c\left(\int_{V}f^{p}(x)\triangle_{V}^{(\gamma+1\rangle p/q-1}(x)dx\right)^{1/p}.
\end{equation}

If $1\leq p<\infty$ and $\alpha<1/p-\sigma(V)/p^{\prime }$, then

\begin{equation}\label{42}
\displaystyle \mathrm{e}\mathrm{s}\mathrm{s}\sup_{x\in
V}(\triangle_{V}^{-1+\alpha}(x)L(f\triangle_{V}^{-\alpha})(x))\leq
c\left(\int_{V}f^{p}(x)\triangle_{V}^{-1}(x)dx\right)^{1/p}.
\end{equation}

If $\alpha<-\sigma(V),$ then

\begin{equation}\label{43}
\displaystyle \mathrm{e}\mathrm{s}\mathrm{s}\sup_{x\in
V}(\triangle_{V}^{-1+\alpha}(x)L(f\triangle_{V}^{-\alpha})(x))\leq c\; \mathrm{e%
}\mathrm{s}\mathrm{s}\sup_{x\in V}f(x).
\end{equation}
\end{theo}

\begin{proof} We will consider the kernel of the Laplace's operator $k(x,\
y)=e^{-x^{*}\cdot y}$. It was shown earlier that $k(x,\ y)$ is $V\times V$
-homogeneous of order $0$. We will apply Theorem 3.3. To verify (\ref{17}) we note
that it is shown in [7] that if $\delta>\sigma(V),$ then 
\begin{equation*}
\int_{V}e^{-x^{*}\cdot y}\triangle_{V}^{\delta}(y)dy<\infty.
\end{equation*}
To verify (\ref{18}) we note 
\begin{equation*}
\int_{V}e^{-(x^{*}\cdot y)q/p}\triangle_{V}^{\gamma-q
+(\delta+1)q/p^{\prime }}(x)dx
=c\int_{V^{*}}e^{-(x^{*}\cdot y)q/p}\Delta_{V^{*}}^{-\gamma+q-(%
\delta+1)q/p^{\prime }-2}(x^{*})dx^{*}
\end{equation*}
We therefore need $\delta$ which satisfies both
$$
\sigma(V)<\delta \text { and }-\gamma+q-(\delta+1) q / p^{\prime}-2>\sigma\left(V^{*}\right).
$$
The condition $\gamma<-\sigma(V) q / p^{\prime}-\sigma\left(V^{*}\right)+q / p-2$ implies the existence of such $\delta$, and (\ref{41}) is proved.

To prove (\ref{42}), we verify that the kernel
\begin{equation}\label{44}
k(x, y)=e^{-x^{*} \cdot y} \Delta_{V}^{\alpha}(x) \Delta_{V}^{-\alpha}(y)
\end{equation}
satisfies (\ref{23}). The kernel $k$ is $V\times V$-homogeneous of order $0$ and
we have that

$$
\int_{V} e^{-\left(x^{*} \cdot y\right) p^{\prime}} \Delta_{V}^{\alpha p^{\prime}}(x) \Delta_{V}^{-\alpha p^{\prime}}(y) \Delta_{V}^{p^{\prime}-1}(y) d y=\Delta_{V}^{\alpha p^{\prime}}(x) \int_{V} e^{-\left(x^{*} \cdot y\right) p^{\prime}} \Delta_{V}^{-\alpha p^{\prime}+p^{\prime}-1}(y) d y
$$
and the last integral is finite since $\alpha<1 / p-\sigma(V) / p^{\prime}$. The result now follows from Theorem $3.5$.

Finally, to prove (\ref{43}), we verify that the kernel (\ref{44}) satisfies the condition (\ref{25}) with $\delta=0$. Thus,
$$
\int_{V} e^{-x^{*} \cdot y} \Delta_{V}^{\alpha}(x) \Delta_{V}^{-\alpha}(y) d y=\Delta_{V}^{\alpha}(x) \int_{V} e^{-x^{*} \cdot y} \Delta_{V}^{-\alpha}(y) d y
$$
The latter integral is finite since $-\alpha>\sigma(V)$ and the result follows from Theorem 3.6.
\end{proof}

\end{document}